\documentclass[11pt,a4paper]{article}
\usepackage{amssymb,amsmath,amscd,amsthm,amsxtra,latexsym}
\theoremstyle{plain}
\newtheorem{theorem}{Theorem}

\newtheorem{proposition}{Proposition}
\newtheorem{lemma}{Lemma}
\newtheorem{corollary}{Corollary}
\theoremstyle{definition}
\newtheorem{definition}{Definition}
\theoremstyle{remark}
\newtheorem{remark}{Remark}

\newtheorem{example}{Example}

\begin{document}
\title{The extrinsic holonomy Lie algebra of a parallel submanifold}
\author{Tillmann Jentsch}
\date{\today}
\maketitle
\sloppy
% RZ-Makros.TEX
% Dies sind Makros, die fŸr alle mathematischen Texte gut sind.

\sloppy
\renewcommand{\thefootnote}{\fnsymbol{footnote}}
\renewcommand{\labelenumi}{(\alph{enumi})}
\renewcommand{\labelenumii}{(\roman{enumii})}
\renewcommand{\labelenumiii}{(\arabic{enumiii})}

%  Allgemeine Makros

\newcommand{\R}{\mathrm{I\!R}}
\newcommand{\N}{\mathrm{I\!N}}
\newcommand{\HH}{\mathrm{I\!H}}
\newcommand{\F}{\mathrm{I\!F}}
\newcommand{\E}{\mathrm{I\!E}}
\newcommand{\K}{\mathrm{I\!K}}
\newcommand{\PP}{\mathrm{I\!P}}
\newcommand{\Z}{\mathbb{Z}}
\newcommand{\Q}{\mathbb{Q}}
\newcommand{\C}{\mathbb{C}}
\newcommand{\OO}{\mathbb{O}}

\newcommand{\scrA}{\mathcal{A}}
\newcommand{\scrB}{\mathcal{B}}
\newcommand{\scrC}{\mathcal{C}}
\newcommand{\scrD}{\mathcal{D}}
\newcommand{\scrE}{\mathcal{E}}
\newcommand{\scrF}{\mathcal{F}}
\newcommand{\scrG}{\mathcal{G}}
\newcommand{\scrH}{\mathcal{H}} % alt \Hscr
\newcommand{\scrI}{\mathcal{I}}
\newcommand{\scrJ}{\mathcal{J}}
\newcommand{\scrK}{\mathcal{K}}
\newcommand{\scrL}{\mathcal{L}}
\newcommand{\scrM}{\mathcal{M}}
\newcommand{\scrN}{\mathcal{N}}
\newcommand{\scrO}{\mathcal{O}}
\newcommand{\scrP}{\mathcal{P}}
\newcommand{\scrQ}{\mathcal{Q}}
\newcommand{\scrR}{\mathcal{R}}
\newcommand{\scrS}{\mathcal{S}}
\newcommand{\scrT}{\mathcal{T}}
\newcommand{\scrU}{\mathcal{U}}
\newcommand{\scrV}{\mathcal{V}} % alt \Vscr
\newcommand{\scrW}{\mathcal{W}}
\newcommand{\scrX}{\mathcal{X}}
\newcommand{\scrY}{\mathcal{Y}}
\newcommand{\scrZ}{\mathcal{Z}}

\newcommand{\frakA}{\mathfrak{A}}
\newcommand{\frakB}{\mathfrak{B}}
\newcommand{\frakC}{\mathfrak{C}}
\newcommand{\frakD}{\mathfrak{D}}
\newcommand{\frakE}{\mathfrak{E}}
\newcommand{\frakF}{\mathfrak{F}}
\newcommand{\frakG}{\mathfrak{G}}
\newcommand{\frakH}{\mathfrak{H}}
\newcommand{\frakI}{\mathfrak{I}}
\newcommand{\frakJ}{\mathfrak{J}}
\newcommand{\frakK}{\mathfrak{K}}
\newcommand{\frakL}{\mathfrak{L}}
\newcommand{\frakM}{\mathfrak{M}}
\newcommand{\frakN}{\mathfrak{N}}
\newcommand{\frakO}{\mathfrak{O}}
\newcommand{\frakP}{\mathfrak{P}}
\newcommand{\frakQ}{\mathfrak{Q}}
\newcommand{\frakR}{\mathfrak{R}}
\newcommand{\frakS}{\mathfrak{S}}
\newcommand{\frakT}{\mathfrak{T}}
\newcommand{\frakU}{\mathfrak{U}}
\newcommand{\frakV}{\mathfrak{V}}
\newcommand{\frakW}{\mathfrak{W}}
\newcommand{\frakX}{\mathfrak{X}}
\newcommand{\frakY}{\mathfrak{Y}}
\newcommand{\frakZ}{\mathfrak{Z}}

\newcommand{\fraka}{\mathfrak{a}}
\newcommand{\frakb}{\mathfrak{b}}
\newcommand{\frakc}{\mathfrak{c}}
\newcommand{\frakd}{\mathfrak{d}}
\newcommand{\frake}{\mathfrak{e}}
\newcommand{\frakf}{\mathfrak{f}}
\newcommand{\frakg}{\mathfrak{g}}
\newcommand{\frakh}{\mathfrak{h}}
\newcommand{\fraki}{\mathfrak{i}}
\newcommand{\frakj}{\mathfrak{j}}
\newcommand{\frakk}{\mathfrak{k}}
\newcommand{\frakl}{\mathfrak{l}}
\newcommand{\frakm}{\mathfrak{m}}
\newcommand{\frakn}{\mathfrak{n}}
\newcommand{\frako}{\mathfrak{o}}
\newcommand{\frakp}{\mathfrak{p}}
\newcommand{\frakq}{\mathfrak{q}}
\newcommand{\frakr}{\mathfrak{r}}
\newcommand{\fraks}{\mathfrak{s}}
\newcommand{\frakt}{\mathfrak{t}}
\newcommand{\fraku}{\mathfrak{u}}
\newcommand{\frakv}{\mathfrak{v}}
\newcommand{\frakw}{\mathfrak{w}}
\newcommand{\frakx}{\mathfrak{x}}
\newcommand{\fraky}{\mathfrak{y}}
\newcommand{\frakz}{\mathfrak{z}}

\newcommand{\rmA}{\mathrm{A}}
\newcommand{\rmB}{\mathrm{B}}
\newcommand{\rmC}{\mathrm{C}}
\newcommand{\rmD}{\mathrm{D}}
\newcommand{\rmE}{\mathrm{E}}
\newcommand{\rmF}{\mathrm{F}}
\newcommand{\rmG}{\mathrm{G}}
\newcommand{\rmH}{\mathrm{H}}
\newcommand{\rmI}{\mathrm{I}}
\newcommand{\rmJ}{\mathrm{J}}
\newcommand{\rmK}{\mathrm{K}}
\newcommand{\rmL}{\mathrm{L}}
\newcommand{\rmM}{\mathrm{M}}
\newcommand{\rmN}{\mathrm{N}}
\newcommand{\rmO}{\mathrm{O}}
\newcommand{\rmP}{\mathrm{P}}
\newcommand{\rmQ}{\mathrm{Q}}
\newcommand{\rmR}{\mathrm{R}}
\newcommand{\rmS}{\mathrm{S}}
\newcommand{\rmT}{\mathrm{T}}
\newcommand{\rmU}{\mathrm{U}}
\newcommand{\rmV}{\mathrm{V}}
\newcommand{\rmW}{\mathrm{W}}
\newcommand{\rmX}{\mathrm{X}}
\newcommand{\rmY}{\mathrm{Y}}
\newcommand{\rmZ}{\mathrm{Z}}

\newcommand{\bbA}{\mathbb{A}}
\newcommand{\bbB}{\mathbb{B}}
\newcommand{\bbC}{\mathbb{C}}
\newcommand{\bbD}{\mathbb{D}}
\newcommand{\bbE}{\mathbb{E}}
\newcommand{\bbF}{\mathbb{F}}
\newcommand{\bbG}{\mathbb{G}}
\newcommand{\bbH}{\mathbb{H}}
\newcommand{\bbI}{\mathbb{I}}
\newcommand{\bbJ}{\mathbb{J}}
\newcommand{\bbK}{\mathbb{K}}
\newcommand{\bbL}{\mathbb{L}}
\newcommand{\bbM}{\mathbb{M}}
\newcommand{\bbN}{\mathbb{N}}
\newcommand{\bbO}{\mathbb{O}}
\newcommand{\bbP}{\mathbb{P}}
\newcommand{\bbQ}{\mathbb{Q}}
\newcommand{\bbR}{\mathbb{R}}
\newcommand{\bbS}{\mathbb{S}}
\newcommand{\bbT}{\mathbb{T}}
\newcommand{\bbU}{\mathbb{U}}
\newcommand{\bbV}{\mathbb{V}}
\newcommand{\bbW}{\mathbb{W}}
\newcommand{\bbX}{\mathbb{X}}
\newcommand{\bbY}{\mathbb{Y}}
\newcommand{\bbZ}{\mathbb{Z}}

%% Es werden Fu§not-Symbolgleichungsnummerierung und Normal-Gleichungsnummerierung
%% eingefŸhrt. Die ZŠhler sind seqcounter und neqcouter. In jedem neuen Abschnitt
%% beginnt die ZŠhlung von vorne
\newcounter{neqcounter}[section]
\newcounter{seqcounter}[section]
\newcommand{\numberequation}[1]
   {\renewcommand{\theequation}{\arabic{equation}}
    \setcounter{equation}{\value{neqcounter}}
    \begin{equation}#1\end{equation}\stepcounter{neqcounter}}
\newcommand{\symbolequation}[1]
   {\renewcommand{\theequation}{\fnsymbol{equation}}
    \setcounter{equation}{\value{seqcounter}}
    \begin{equation}#1\end{equation}\stepcounter{seqcounter}}

\newcommand{\diff}{\mathrm{d}}
\newcommand{\Diff}{\mathrm{D}}
\newcommand{\id}{\mathrm{id}}
\newcommand{\GL}{\mathrm{GL}}
\newcommand{\SL}{\mathrm{SL}}
\newcommand{\SO}{\mathrm{SO}}
\newcommand{\End}{\mathrm{End}}
\newcommand{\Alt}{\mathrm{Alt}}
\newcommand{\Spur}{\mathrm{Spur}}
\newcommand{\rk}{\mathrm{rk}}
\newcommand{\Kern}{\mathrm{Kern}}
\newcommand{\Spann}{\mathrm{Spann}}
\newcommand{\Top}{\mathrm{Top}}
\newcommand{\del}{\partial}
\newcommand{\ddel}[2]{\frac{\partial#1}{\partial#2}}
\newcommand{\ddeldel}[3]{\frac{\partial^2#1}{\partial#2\:\partial#3}}
\newcommand{\ddelsquare}[2]{\frac{\partial^2#1}{\partial{#2}^2}}
\newcommand{\grad}{\mathrm{grad}}
\newcommand{\hess}{\mathrm{hess}}
\newcommand{\Hess}{\mathrm{Hess}}
\newcommand{\Rp}{\R_+}
\newcommand{\eps}{\varepsilon}
\newcommand{\vi}{\varphi}
\newcommand{\vkap}{\varkappa} %alt
\newcommand{\kap}{\varkappa}  %neu
\newcommand{\thet}{\vartheta}
\newcommand{\ro}{\varrho}
\newcommand{\Hol}{\mathrm{Hol}}
\newcommand{\rg}{\mathrm{rg}}

\renewcommand{\O}{\varnothing}

\newcommand{\UL}{\mathchoice
{\setbox0=\hbox{$\displaystyle U$}
 \setbox1=\hbox{$\displaystyle l$}
 \hbox{\box0\kern-0.80\wd1\lower0.029\ht1\box1}}
{\setbox0=\hbox{$\displaystyle U$}
 \setbox1=\hbox{$\displaystyle l$}
 \hbox{\box0\kern-0.80\wd1\lower0.029\ht1\box1}}
{\setbox0=\hbox{$\scriptstyle U$}
 \setbox1=\hbox{$\scriptstyle l$}
 \hbox{\box0\kern-0.80\wd1\lower0.029\ht1\box1}}
{\setbox0=\hbox{$\scriptscriptstyle U$}
 \setbox1=\hbox{$\scriptscriptstyle l$}
 \hbox{\box0\kern-0.80\wd1\lower0.029\ht1\box1}}
}

\newcommand{\Uo}{\UL^o}
\def\fam#1#2#3#4{({#1}_{#2})_{#2=#3,\dotsc,#4}}
\def\zz#1#2#3{#1=#2,\dotsc,#3}
\newcommand{\intint}[2]{\{#1,\dotsc,#2\}}
\newcommand{\Vektor}[2]{({#1}_1,\dotsc,{#1}_{#2})} 

\newcommand{\qmq}[1]{\quad\mbox{#1}\quad}
\newcommand{\Menge}[2]{\{\,#1\,|\,#2\,\}}
\newcommand{\Abstand}[1]{\mbox{}\par\vspace{-#1mm}}
\newcommand{\ninfty}{{-}\infty}

\def\bild#1#2#3#4{\leavevmode\vbox to #1{\vfil \hbox to #2{\special{picture #4
   scaled #3}\hfil}}}
   %Parametr: #1=Hšhe,#2=Breite,#3=Skalierung 0.1<f<10,#4=Bildname

%%  Zu den Rahmen:
%%  grRahmen ist mit rahmen aus dem Alalysisscript identisch.
%%  grRahmen geht fast Ÿber die ganze Textbreite
%%
%%  klRahmen ist aus kleinerRahmen aus den RZ-Makros von 1996 entstanden;
%%  seine Hšhe ist etwas vergrš§ert worden,
%%  seine Breite ist dem Text angepasst, dieser darf nur eine Zeile sein
%%
%%  varRahmen ist eine Neuschšpfung: Man stellt eine Breite ein. Dann wird der
%%  in dem Rahmen wie in einer Minipage behandelt.
%%
\newcommand{\normRahmen}[1]
           {\setlength{\fboxsep}{8pt}\fbox{#1}}

\newcommand{\grRahmen}[1]{\begin{center}\setlength{\fboxrule}{0.8pt}
           \setlength{\fboxsep}{8pt}
              \fbox{\begin{minipage}{140mm}\rule{0mm}{5mm}\hspace*{-2mm} #1
                    \end{minipage}}\end{center}}

\newcommand{\klRahmen}[1]{\begin{center}\large\setlength{\fboxrule}{0.8pt}
                     \setlength{\fboxsep}{8pt}
                     \fbox{\rule[-2mm]{0mm}{7mm} #1 }
                     \normalsize\end{center}}

\newcommand{\varRahmen}[2]{\begin{center}\setlength{\fboxrule}{0.8pt}
           \setlength{\fboxsep}{8pt}
              \fbox{\begin{minipage}{#1}\rule{0mm}{5mm}\hspace*{-2mm} #2
                    \end{minipage}}\end{center}}

% die 1996-Version
%\newcommand{\kleinerRahmen}[1]{\begin{center}\large\setlength{\fboxrule}{0.8pt}
%                     \setlength{\fboxsep}{8pt}\fbox{#1}\normalsize\end{center}}

%% Randnotizen vom Juli 1997
\setlength{\marginparsep}{5pt}
\setlength{\marginparwidth}{30pt}
\newcommand{\marginlabel}[1]           
                 {\mbox{}\marginpar{\raggedleft\hspace{0pt}#1}}
\newcommand{\Ausrufezeichen}{\marginlabel{\raisebox{-1.2ex}{\Huge$\boldsymbol{!}$}}}
\newcommand{\Randbalken}[2]{\marginlabel{{\rule[#1]{0.5mm}{#2}}}}
\newcommand{\Zeigefinger}{\marginlabel{\raisebox{-0.8ex}{\LARGE\ding{43}}}}
\newcommand{\Blume}{\marginlabel{\raisebox{-0.6ex}{\LARGE\ding{94}}}}

%  Differentialgeometrie - Makros
\newcommand{\bigoperp}{\mathop{\bigcirc\raisebox{0.25em}
 {\hskip-0.53em\hbox{\vrule height1.0ex width0.04em}
  \hskip-0.77em\hbox{\vrule height0.04em width 0.8em}
  \hskip- 0.2em}}}

\renewcommand{\bigoplus}{\mathop{\bigcirc
  \raisebox{-0.22em}{\hskip-0.53em\hbox{\vrule height2.08ex width0.04em}
  \raisebox{ 0.48em}{\hskip-0.75em\hbox{\vrule height0.04em width 0.8em}}
  \hskip- 0.2em}}}

\def\NP#1#2#3{\mathchoice
  {{\textstyle{\prod\limits_{i=#2}^{#3}}{#1}_i}} 
  {{\textstyle{\prod_{i=#2}^{#3}}{#1}_i}}
  {}
  {}
   }
\def\DS#1#2#3{\bigoplus_{i=#2}^{#3}\!#1_i} 
\def\OS#1#2#3{\bigoperp_{i=#2}^{#3}#1_i} 
\def\WP#1#2#3{#1_0 \times_{#2}\NP {#1}{1}{#3}}
\def\TP#1#2#3#4{{\rule{0mm}{2ex}}^{#2}\NP{#1}{#3}{#4}}

\newcommand {\cinf}{\ensuremath{\mathrm{C}^{\infty}}}
\newcommand {\X}{\mathfrak{X}}
\newcommand{\Tensor}[3]{\mathfrak{T}^{(#2,#3)}(#1)}
\newcommand{\alphad}{\dot{\alpha}}
\newcommand {\g}[2]{\langle #1,#2\rangle}
\def\gg{\g{\cdot\,}{\cdot}}
\newcommand {\euc}[2]{\langle\!\langle #1,#2 \rangle\!\rangle}
\newcommand {\peuc}[2]{\langle\!\langle #1,#2 \rangle\!\rangle\!_s}
\newcommand {\Kov}[2]{\nabla_{#1}#2}
\newcommand {\Kovh}[2]{\widehat{\nabla}_{#1}#2}
\newcommand {\Kovt}[2]{\widetilde{\nabla}_{#1}#2}
\newcommand {\Kovperp}[2]{\nabla^{\perp}_{#1}#2}
\newcommand {\Kovv}[3]{{\nabla^{\scriptscriptstyle{#1}}}_{\!#2}#3}
\newcommand{\Kodel}[1]{\Kov{\del}{#1}}
\newcommand{\Kokodel}[1]{\Kov{\del}{\Kov{\del}{#1}}}

\newcommand {\operp}{\mathbin{\mbox{$\ominus\raisebox{2.9pt}
 {\hskip-0.42em\hbox{\vrule height0.7ex width0.02em}\hskip0.42em }$}}}
\newcommand {\TpM}{T_{p}M}
\newcommand {\Tpf}{T_{p}f}
\newcommand {\Nf}{\perp\!\!(f)}
\newcommand {\NM}{\perp\!\!M}
\newcommand {\Npf}{\perp_p\!\!(f)}
\newcommand {\NpM}{\perp_p\!\!M}
\newcommand {\Nepf}{\perp^{\!\!1}_p\!\!(f)}
\newcommand {\NepM}{\perp^{\!\!1}_p\!\!M}
\newcommand{\Shop}[3]{{\mathrm{A}^{^{\scriptscriptstyle{\!\!#1}}}_{#2}}#3} %shape op

\def\tscrV{\widetilde{\mathcal{V}}}
\def\tscrH{\widetilde{\mathcal{H}}}
\newcommand{\hdisp}[3]{\overset{#2}{\underset{#1}{\parallel}}\!\!#3\,} 
    % horicontal displacement

% Tangentialvektoren
\newcommand{\dotdiff}[2]{\dot{\frac{\diff\ }{\diff #1}}\Big|_{#1=0}\big(#2\big)}
    % NŸbel-Bezeichnung fŸr Tangentialvektoren an eine Kurve t -> a(t)
\newcommand{\dotdif}[2]{\dot{\frac{\diff\ }{\diff #1}}\Big|_{#1=0}#2}
    % NŸbel-Bezeichnung fŸr Tangentialvektoren an eine Kurve a(t)

% Entsprechendes fŸr normale Ableitungen
\newcommand{\normaldiff}[3]{\frac{\diff\ }{\diff #1}\Big|_{#1=#2}\big(#3\big)}
    % NŸbel-Bezeichnung fŸr Ableitungen einer Funktion t -> f(t)
\newcommand{\normaldif}[3]{\frac{\diff\ }{\diff #1}\Big|_{#1=#2}#3}
    % NŸbel-Bezeichnung fŸr Ableitungen von f(t)

% Reelle projektive RŠume
\newcommand {\RPn}{\ensuremath{\R\mathrm{P}^{n}}}
\newcommand {\RPm}{\ensuremath{\R\mathrm{P}^{m}}}

% Komplexe Raumformen
\newcommand {\CPn}{\ensuremath{\C\mathrm{P}^{n}}}
\newcommand {\CPm}{\ensuremath{\C\mathrm{P}^{m}}}
\newcommand {\CHn}{\ensuremath{\C\mathrm{H}^{n}}}
\newcommand {\CHm}{\ensuremath{\C\mathrm{H}^{m}}}

\hyphenation{Riemann-sche Riemann-ian}

%DiffGeo-Makros 2002
% Hier werden Markos zusammengestellt, die fŸr das Diffgeo-Skript nŸtzlich sind

\newcommand{\sst}{\scriptscriptstyle}
\newcommand{\VV}{\mathbb{V}}

% Makros fŸr DiffGeo-AnhE

\newcommand{\PV}{\rmP(\bbV)}
\newcommand{\SV}{\rmS(\bbV)}
\newcommand{\EPV}{\E\times\PV}
\newcommand{\ESV}{\E\times\SV}

% SpezialitŠten

\newcommand{\GrTM}{\rmG_r(TM)}  % Grassmann-BŸndel Ÿber M
\newcommand{\Kovp}[1]{\Kov{}{#1}\big|_{p_0}} 
    % interessant fŸr affine Vektorfelder, Skript Ÿber homogene RŠume, Karcher-Theorie

\newcommand{\GrnTM}{\rmG_r^\nabla(TM)}
\newcommand{\MtM}{{M^\times}}

%Till-Makros 2004
%

\newcommand{\ghdisp}[4]{(\hdisp{#1}{#2}{#3})^{\sst{#4}}}
\newcommand{\hh}[1]{\hat{\!\hat{#1}}} %% <--
\renewcommand{\R}{\bbR}
\renewcommand{\Kodel}[2]{\nabla^
  {\text{\raisebox{0.5ex}{${\scriptstyle{#1}}$}}}_{\del}{#2}}
\renewcommand{\Kokodel}[2]{\Kodel{#1}{\Kodel{#1}{#2}}}
\renewcommand{\E}{\bbE}
\renewcommand{\Kovh}[2]{\hat\nabla_{#1}{#2}}
\newcommand{\WBKov}[2]{\overline{\nabla}_{#1}{#2}}
\newcommand{\Till}{\ensuremath {\scrT}}
\newcommand{\Tillb}{\ensuremath {\Till_b}}
\newcommand{\Tsukada}{\ensuremath {\scrD}}
\newcommand{\Tsukadaq}{\ensuremath {\Tsukada_q}}
\newcommand{\Gzwei}{\ensuremath {\rmG_m^{2}(N)}}
\newcommand{\Isom}{(\hat\eta,\hat\nu)}
\newcommand{\hv}{\hat{v}}
\newcommand{\hu}{\hat{u}}
\newcommand{\hN}{\hat{N}}
\newcommand{\hnu}{\hat\nu}
\newcommand{\hW}{\hat{W}}
\newcommand{\hX}{\hat{X}}
\newcommand{\hY}{\hat{Y}}
\newcommand{\hx}{\hat{x}}
\newcommand{\hy}{\hat{y}}
\newcommand{\hU}{\hat{U}}
\newcommand{\hV}{\hat{V}}
\newcommand{\hZ}{\hat{Z}}
\newcommand{\hOmega}{{\hat{\Omega}}}

\newcommand{\GamA}{\Gamma_{\!\!A}}
\newcommand{\piP}{\pi^{\sst\bbP}}
\newcommand{\rmFB}{\mathrm{FB}}
\newcommand{\Besym}{\scrB_{esym}(N)}
\newcommand{\Besymo}{\scrB_{esym}^0(N)}
\newcommand{\ttop}{{\boldsymbol{\top\hspace{-0.75em}\top}}}
\newcommand{\bbot}{{\boldsymbol{\!\bot\hspace{-0.75em}\bot}}}

\newcommand{\kernel}{\mathrm{kernel}}

\newcommand{\Gl}{\mathrm{Gl}}
\newcommand{\SU}{\mathrm{SU}}
\newcommand{\Sp}{\mathrm{Sp}}
\newcommand{\so}{\mathfrak{so}}
\newcommand{\su}{\mathfrak{su}}
\newcommand{\gl}{\mathfrak{gl}}
\renewcommand{\sp}{\mathfrak{sp}}
\newcommand{\tr}{\mathfrak{tr}}
\newcommand{\osc}{\scrO}
\newcommand{\Id}{\mathrm{Id}}
\newcommand{\Iso}{\mathrm{I}}
\newcommand{\trace}{\mathrm{trace}}
\newcommand{\Sym}{\mathrm{Sym}}
\renewcommand{\i}{\mathrm{i}}
\renewcommand{\j}{\mathrm{j}}
\renewcommand{\k}{\mathrm{k}}
\renewcommand{\Spann}[2]{\{#1\big|#2\}_{\scriptstyle\R} }\,
\newcommand{\killing}{\mathfrak{f}}
\newcommand{\hol}{\mathfrak{hol}}
\newcommand{\fetth}{\boldsymbol{h}}
\newcommand{\fettb}{\boldsymbol{b}}
\newcommand{\ad}{\mathrm{ad}}
\newcommand{\Ad}{\mathrm{Ad}}
\newcommand{\Hom}{\mathrm{Hom}}
\newcommand{\Spec}{\mathrm{Spec}}
\newcommand{\codim}{\mathrm{codim}}
\newcommand{\rank}{\mathrm{rank}}

\abstract{We investigate parallel submanifolds of a Riemannian symmetric
space $N$\,. The special case of a symmetric submanifold has been investigated
by many authors before and is well understood.
We observe that there is an intrinsic property of the second fundamental form which
distinguishes full symmetric submanifolds from arbitrary full parallel
submanifolds of $N$\,, usually called ``1-fullness of $M$''\,. Furthermore,
for every parallel submanifold $M\subset N$ we consider the pullback bundle
$TN|M$ with its induced connection, 
which admits a distinguished parallel subbundle
$\osc M$\,, usually called the ``second osculating bundle of $M$''\,. If $M$
is a complete parallel submanifold of $N$\,, then we can describe the
corresponding holonomy Lie algebra of $\osc M$ by means of the second 
fundamental form of $M$ and the curvature tensor of $N$ at the origin\,.
If moreover $N$ is simply connected and $M$ is even a full symmetric submanifold 
of $N$\,, then we will calculate the holonomy Lie algebra of $TN|M$ in an explicit form.}

% \cleardoublepage
%\tableofcontents
\section{Introduction}
In this article, $N$ denotes a Riemannian symmetric space. For an isometric
immersion $f:M\to N$\,, let $TM$\,, $\bot f$\,, $h:TM\times TM\to\bot f$ and $S:TM\times\bot f\to TM$ denote the tangent bundle
of $M$\,, the normal bundle of $f$\,, the second
fundamental form and the shape operator, respectively. Let $\nabla^M$ and $\nabla^N$ denote the
Levi Civita connection of $M$ resp.\ of $N$ and $\nabla^\bot$ the usual
connection on $\bot f$ (obtained by projection).
The equations of Gau{\ss} and Weingarten state for $X,Y\in \Gamma(TM), \xi\in\Gamma(\bot f)$
\begin{equation}\label{eq:GW}
\Kovv{N}{X}{Tf\,Y}=Tf(\Kovv{M}{X}{Y})+h(X,Y)\qmq{and}\Kovv{N}{X}{\xi}=-Tf(S_\xi(X))+\Kovv{\bot}{X}{\xi}\;.
\end{equation}
On the vector bundle $\rmL^2(TM,\bot f)$ there is a connection induced by
$\nabla^M$ and $\nabla^\bot$ in a natural way, often called ``Van der Waerden-Bortolotti connection''.

\bigskip
\begin{definition}\label{de:parallel} $f$ is called
 \emph{parallel} if its second
fundamental form $h$ is a parallel section of
the vector bundle $\rmL^2(TM,\bot M)$\,.
\end{definition}
In a similar fashion, we define parallel submanifolds of $N$ (via the isometric immersion
given by the inclusion map $\iota^M:M\hookrightarrow N$). 

\bigskip
\begin{example}[Circles]\label{ex:circles}
A unit speed curve $c:J\to N$ is parallel if and only if it satisfies the equation
\begin{equation}
    \label{eq:circles}
\nabla^N_\partial\nabla^N_\partial\dot c=-\kappa^2 \dot c
\end{equation}
for some constant $\kappa\in\R$\,. For $\kappa=0$ these curves are geodesics;
otherwise,
due to Nomizu and Yafo in~\cite{NY}, $c$ is called an (extrinsic) circle.
One can show that for every pair $(u,v)\in T_pN\times T_pN$ with $\|u\|=1$ there
exists a unique solution $c$ of~\eqref{eq:circles} defined on the whole real
line with $\dot c(0)=u,\;\Kodel{N}{\dot c}(0)=v$\,.
It is obtained as the envelopment of some straight line or some circle in $T_pN$\,, see also~\cite{JR}\;.
\end{example}

So far, a classification of parallel isometric immersions has been achieved only
if the ambient space is a rank-1 symmetric space.
(see~\cite{BCO},~Ch.\,9.3).
Nevertheless, even if $N$ is of higher rank, then the special case of a 
{\em symmetric} \/submanifold is completely understood  by the work of H.~Naitoh and
others (for an overview on the classification of symmetric submanifolds of
symmetric spaces see~\cite{BCO},~Ch.\,9.4)\,.

\bigskip
\begin{definition}\label{de:extrinsically_symmetric}
$M$ is called a
{\em symmetric submanifold}\/ of $N$ if $M$ is a symmetric space (whose
geodesic symmetries are denoted by $\sigma^M_p\;(p\in M)$) and for every point $p\in M$ there
exists an involutive isometry $\sigma^\bot_p$ of $N$ such that
\begin{itemize}
\item $\sigma^\bot_p(M)=M$\,,
\item
 $\sigma^\bot_p|M=\sigma_p^M$\,,
\item
and the differential $T_p\sigma^\bot_p$ is the linear reflection
 in the normal space $\bot_pM$\,.
\end{itemize}
Then we also say that $M$
is extrinsically symmetric in $N$\,.
The family $\sigma^\bot_p\;(p\in M)$ is unique (if it exists) and
is called the \emph{extrinsic symmetries} of $M$\,.
\end{definition}
In fact, symmetric
submanifolds of $N$ are parallel\,, but the converse is not true. So far there seems to be not much known about
arbitrary (i.e.\ not necessarily extrinsically symmetric) parallel submanifolds of an
irreducible symmetric space $N$ of higher rank,
except for a result of K.~Tsukada~\cite{Ts} on parallel
K{\"a}hler submanifolds of Hermitian symmetric spaces (which, in case $N$ is of
higher rank, can be interpreted as a negative result) and the analogue in~\cite{ADM}
for parallel submanifolds of K{\"a}hlerian type in a quaternionic-K{\"a}hler symmetric space of
non-vanishing scalar curvature. 
% As can be seen from the results presented in~\cite{BCO}, if $N$ is irreducible, then there are not
% known yet any other examples of higher dimensional, full\footnote{A  submanifold $M\subset N$ is called full if it is not contained in any proper, totally geodesic
% submanifold of $N$\,.}, complete parallel submanifolds of $N$ than the symmetric submanifolds of $N$\,.

The aim of this article is three-fold: 
\begin{itemize}
\item
First, we will relate the extrinsic symmetry of a
full parallel submanifold of $N$ to an {\em intrinsic} \/property of the second fundamental
form called ``1-fullness of $M$'' (see Definition~\ref{de:1-full} and Theorem~\ref{th:result_1}).

\item
Second, for every complete
parallel submanifold $M\subset N$ we will introduce the extrinsic holonomy Lie
algebra of $M$ resp.\ of its second osculating bundle (see Definition~\ref{de:extrinsic_holonomy}), and we will 
be able to express the latter Lie algebra only in terms of the second 
fundamental form of $M$ and the curvature tensor of $N$ at the origin (see Theorem~\ref{th:hol}).
% If $M$ is even a full submanifold of $N$\,, then this Lie algebra is
% isomorphic to the extrinsic holonomy Lie algebra of $M$ 

\item
Third, for the full symmetric submanifolds of the simply connected
symmetric spaces we will calculate their extrinsic holonomy Lie algebras in
an explicit form (up to certain exceptions, see Theorem~\ref{th:hol_for_symmetric_M}).
\end{itemize}
The precise definitions and the statement of the theorems can be found in the next
Section. 

\bigskip
In a forthcoming paper~\cite{J1}, the extrinsic
homogeneity of (arbitrary) parallel submanifolds in an ambient
symmetric space of possibly higher rank will be studied, for which
Theorem~\ref{th:hol} of this article will serve as a useful tool; moreover, it seems possible that the explicit calculations in the
extrinsically symmetric case could also be helpful for the further study of arbitrary parallel submanifolds in symmetric spaces.

\bigskip
This article was written at the Mathematical Institute of the University of
Cologne. I would like to thank everybody who supported me in the making of this paper.
Special thanks goes to my teacher Professor H.~Reckziegel for his helpful
advises, which served as a sort of ``nutrient medium'' for this article.

\subsection{Overview}
\label{s:overview}
This section gives a detailed, self contained overview on the results
presented in this article, the necessary notation included. 
In Section~\ref{s:results} we recall some
well known properties of parallel submanifolds, and we consider certain
relevant examples. Given an isometric immersion $f:M\to N$\,, in order to keep our
notation as simple as possible, here and in the following we implicitly
identify the tangent space $T_pM$ with the ``first osculating space''
$Tf(T_pM)$ by means of the injective linear map $T_pf$ for each $p\in
M$\,. Then we introduce for each $p\in M$ 
the {\em first normal
space}
\begin{equation}\label{eq:bot_eins}
\bot^1_pf:=\Spann{h(x,y)}{x,y\in T_pM}\;,
\end{equation}
and the {\em second osculating space}
\begin{equation}\label{eq:osc_M}
\osc_pf:=T_pM\oplus\bot^1_pf\;,
\end{equation}
seen as a linear subspace of $T_{f(p)}N$\,.
If $M\subset N$ is actually a (smoothly embedded) submanifold, then the first normal space
$\bot^1_pM$ and the second osculating space $\osc_pM$ are defined as before
via the isometric immersion $\iota^M:M\hookrightarrow N$\,. 

\bigskip
\begin{definition}\label{de:1-full}
\begin{enumerate}
\item In accordance with~\cite{BCO},~Ch.\,2.5, an isometric immersion $f:M\to N$ is called \emph{full} if $f(M)$ is not contained
in any proper, totally geodesic submanifold $\bar N\subset N$\,.
\item In accordance with~\cite{Ts}, an isometric immersion $f:M\to N$ is called \emph{1-full} if always
the first normal space $\bot^1_pf$ coincides with the normal space $\bot_pf$\,.
\end{enumerate}
\end{definition}
Note that there always exists a smallest complete, totally geodesically submanifold $\bar N\subset N$ which contains
$f(M)$\,, and then $\bar N$ is a symmetric space and $f:M\to\bar N$ is a
full isometric immersion. However, 1-fullness is a somehow more
intrinsic property of $f$\,.

As a consequence of the Gau{\ss}
Equation, we see that 1-fullness implies fullness, but the converse is
not true even for parallel isometric immersions:

\bigskip
\begin{definition} 
Let $R^N$ denote the curvature tensor of $TN$\,. A linear subspace $V\subset T_pN$ is curvature invariant if
$R^N(V,V)\,V\subset V$\,. 
\end{definition}

\bigskip
\begin{example}
\label{ex:CP2}
There exists a full circle $c:\R\to\C\rmP^2$ (see Example~\ref{ex:circles}) which
is not 1-full and whose normal spaces are not curvature invariant.
\end{example}
\begin{proof} 
For $p:=(1:0:0)\in \C\rmP^2$ let $u,v\in T_p\C\rmP^2$ be two vectors with $\|u\|=1$ and the property that
$\{u,v\}_\R$ is neither a totally real nor a complex linear space.
Then there exists a circle $c:\R\to N$ with the initial conditions
$$\dot c(0)=u\qmq{and}\nabla^N_\partial \dot c(0)=v\;.$$
Suppose that $\bar N\subset \C\rmP^2$ is a totally geodesic submanifold such
that $c(\R)\subset\bar N$\,.
Thus $T_p\bar N$ is a curvature invariant linear subspace with $\{u,v\}\subset T_p\bar N$\,. Since all curvature invariant
linear subspaces of $\C\rmP^2$ are either totally real or complex, it follows
by construction that $T_p\bar N=T_p\C\rmP^2$ ; thus $c$ is full. The last statement follows, because the normal spaces of $c$ are
three-dimensional (and hence neither totally real nor complex subspaces).
\end{proof}

Section~\ref{s:osc} deals with the proof of the following theorem, which can
not be found in the literature so far\footnote{When I talked about my results at Augsburg, I learned
  that the result described in
  Part~(a) of Theorem~\ref{th:result_1} could also be found in an unpublished paper by E. Heintze.}:

\bigskip
\begin{theorem}\label{th:result_1}
\begin{enumerate}
\item The first normal
spaces $\bot_p^1f$ of a parallel isometric immersion are always curvature invariant.
\item Let a simply connected symmetric space $N$ and a submanifold $M\subset N$
be given. $M$ is a full symmetric submanifold of $N$ if and only if $M$ is a
1-full, complete\footnote{According to
Theorem~7 of~\cite{JR}, for every (not
necessarily complete) parallel submanifold $M_{loc}\subset N$ there exists a simply connected
Riemannian symmetric space $M$\,,
a parallel isometric immersion $f:M\to N$ and an open subset $U\subset M$\,,
such that $f|U:U\to M_{loc}$ is covering. Hence, loosely
said, all parallel submanifolds
can be ``extended'' to simply connected, complete, immersed parallel submanifolds and therefore the completeness assumption in the above theorem
is not too striking.}, parallel submanifold of $N$\,. 
\item Let $N$ be a simply connected symmetric space, which has no Euclidian
  factor (in the sense of the ``de Rham decomposition theorem'',
  see~\cite{BCO},~p.\,290). If
$M$ is a full symmetric submanifold of $N$\,, then at each point $p\in M$
the second fundamental form $h_p$ is a non-degenerate symmetric bilinear form. 
\end{enumerate}
\end{theorem}

If $V$ is a curvature invariant subspace of $T_pN$\,, then $\exp^N(V)\subset N$ (where $\exp^N$ denotes the
exponential spray defined on $TM$) is a totally geodesic
submanifold by a result due to E.~Cartan.
The following result on the ``reduction of the
codimension'' (in the sense of~\cite{Er}) is well known (cf. Lemma 2.1 of~\cite{Ts}); 
it is in fact a consequence of Theorem 3.4 in~\cite{D}
combined with Part~(d) of Proposition~\ref{p:properties} further below:

\bigskip
\begin{theorem}[Dombrowski]\label{th:Dombi}
If $f:M\to N$ is parallel and if at one point $p\in M$ the second osculating space
$\osc_pf$ is contained in some curvature invariant subspace $V\subset T_pN$\,, then
$f(M)\subset\bar N$\,, where $\bar N$ denotes the totally geodesic submanifold $\exp_p(V)\subset N$\, (which again is a symmetric space).
\end{theorem}
Combining Theorem~\ref{th:result_1} with Theorem~\ref{th:Dombi} we hence obtain:

\bigskip
\begin{corollary}\label{co:reduction}
If at one point $p\in M$ the linear space $\osc_pM$ of a parallel submanifold
is a curvature invariant
subspace of $T_pN$ and if then $\bar N:=\exp^N_p(\osc_pM)$ is simply connected,
then $M$ is an extrinsically symmetric submanifold of $\bar N$\,. 
\end{corollary}

Corollary~\ref{co:reduction} should be compared with Lemma 2.2 of~\cite{Ts}.
But note that the second osculating spaces of a parallel isometric immersion
are not {\em not}\ always curvature-invariant (see Example~\ref{ex:CP2}); hence Corollary~\ref{co:reduction} is not always applicable. 

\bigskip
Parallel submanifolds are sometimes also called ``weakly locally symmetric
submanifolds'' (cf.~\cite{NT}). Theorem~\ref{th:wes} in
Section~\ref{sec:symmetric_strip} will give a geometric reason for that notion, as follows:

For every parallel isometric immersion $f:M\to N$ 
we introduce the pullback bundle $$f^*TN:=\underset{p\in
  M}{\bigcup}\{p\}\times T_{f(p)}N$$ (which is a vector bundle over $M$);
moreover, $\nabla^N$ defines a connection on $f^*TN$\,.
According to Proposition~\ref{p:properties}, the {\em second osculating bundle} $$
\osc f:=\underset{p\in M}{\bigcup}\{p\}\times\osc_pf$$ is a $\nabla^N$-parallel subbundle of
$f^*TN$\,, hence $\osc f$ is equipped with the connection $\nabla^{\osc f}$
induced by restriction of $\nabla^N$\,. If $f:M\to N$ is a parallel isometric immersion defined on a simply
connected symmetric space $M$ (cf. Proposition~\ref{p:locally_symmetric})\,,
then we can proof the existence of certain distinguished vector bundle involutions on
$\osc f$\,; in this way, we finally come to the conclusion that $M$ is ``extrinsically symmetric in $\osc f$'' (in a weak sense). 
However, due to its technical nature the precise statement of 
Theorem~\ref{th:wes} is skipped at this point of the paper.
As a first consequence of Theorem~\ref{th:wes}, we will see that $\bot^1f$ is a
{\em homogeneous vector bundle}\/ over $M$ (Proposition~\ref{p:homogeneity_of_the_first_normal_bundle}).

\bigskip
We now introduce the extrinsic holonomy Lie algebras of a parallel isometric
immersion $f:M\to N$ and of its second osculating bundle with respect to some base
point $o\in M$\,. 
For each differentiable curve $c:[0,1]\to N$ let
$\ghdisp{0}{1}{c}{N}$ denote the parallel displacement in $TN$ along $c$ and consider the Holonomy groups of
$TN$ with respect to $\nabla^N$ and of $\osc f$ with respect to $\nabla^{\osc f}$ (the connection
which was introduced above), respectively:
\begin{align}
\label{eq:Hol(N)}
&\Hol(N):=\big\{\ghdisp{0}{1}{c}{N}\big{|}c:[0,1]\to N\text{\ is a loop
    with }c(0)=f(o)\big\}\;,\\
\label{eq:Hol(TN|M)}
&\Hol(f^*TN):=\big\{\ghdisp{0}{1}{f\circ c}{N}\big{|}c:[0,1]\to M\text{\ is a loop
    with }c(0)=o\big\}\;,\\
\label{eq:Hol(osc_M)}
&\Hol(\osc f):=\big\{\ghdisp{0}{1}{f\circ c}{\osc f}\big{|}c:[0,1]\to M\text{\ is a loop
    with }c(0)=o\big\}\;.
\end{align}
Then $\Hol(N)$ and $\Hol(f^*TN)$ are known to be Lie subgroups of
$\SO(T_{f(o)}N)$\,, and  $\Hol(\osc f)$ is a Lie subgroup of $\SO(\osc_of)$\,; the corresponding
Lie algebras are denoted by $\hol(N)$ resp.\ by $\hol(f^*TN)$\,. Moreover, $\Hol(f^*TN)\subset\Hol(N)$ is a Lie subgroup and hence
$\hol(f^*TN)$ is a Lie subalgebra of $\hol(N)$\,. If $M\subset N$ is a parallel
submanifold, then we define the pullback bundle $TN|M$ and the second
osculating bundle of $M$ via the isometric immersion
$f=\iota^M$\,. Then the Lie groups $\Hol(TN|M)$ and $\Hol(\osc M)$ and
their Lie algebras $\hol(TN|M)$ and $\hol(\osc M)$ are defined in a similar fashion.

\bigskip
\begin{definition}\label{de:extrinsic_holonomy}
\begin{enumerate}
\item
We will call $\hol(f^*TN)$ resp.\ $\hol(TN|M)$ the {\em extrinsic holonomy Lie algebra}\ of the
immersion $f$ resp.\ of the submanifold $M$\,.
\item $\hol(\osc f)$ resp.\ $\hol(\osc M)$ will be called the {\em extrinsic holonomy Lie algebra}\ of $\osc f$
  resp.\ of $\osc M$\,.
\end{enumerate}
\end{definition}

\bigskip
\begin{example}
Let $M$ be a totally geodesic submanifold of $N$\,. Since both the
vector subbundle $TM\subset TN|M$ and the curvature tensor $R^N$ are parallel with respect to $\nabla^N$\,, the
Theorem of Ambrose/Singer implies that $\hol(TN|M)=\Spann{R^N(x,y)}{x,y\in T_oM}$\,.
\end{example}

\bigskip
\begin{remark}\label{re:hol(TN|M)}
By means of the Theorem of Ambrose/Singer, 
a (parallel) isometric immersion $f:M\to N$ is curvature isotropic 
(i.e. $R^N(x,y)=0$ for all $x,y\in
T_pM$ and $p\in M$\,, cf.~\cite{FP}) if and only if $\hol(f^*TN)=\{0\}$\,. Therefore,
briefly said, $\hol(f^*TN)$ measures ``how far $f$ is away from being curvature
isotropic''\,.
\end{remark}

\bigskip
The next theorem describes the general structure
of $\hol(\osc f)$ only in terms of the curvature tensor $R^N$ at $f(o)$ and the second
fundamental form of $f$ at $o$\,. For this we will need the following
notation: 

For an arbitrary Euclidian vector space $V$ and some subspace $W\subset V$ let $\sigma^\bot\in \rmO(V)$ denote 
the linear reflection in $W^\bot$ and $\Ad(\sigma^\bot):\so(V)\to\so(V),A\mapsto \sigma^\bot\circ A\circ\sigma^\bot$
the induced involution on $\so(V)$\,. 
Let $\so(V)_+$ resp.\ $\so(V)_-$ be the $+1$- resp.\ $-1$-eigenspaces of
$\Ad(\sigma^\bot)$\,, i.e.
\begin{align}
    \label{eq:even}
&\so(V)_+:=\left\{\left.\left(
\begin{array}{cc}
A & 0 \\
0 & B
\end{array}
\right )\right|A\in \so(W),B\in \so(W^\bot)\right\}\;, \\
\label{eq:odd}
&\so(V)_-:=\left\{\left.\left(
\begin{array}{cc}
0 & -C^* \\
C & 0
\end{array}
\right)\right|C\in \rmL(W,W^\bot) \right\}\;.
\end{align}
Then the rules for $\Z /2\Z$ graded Lie algebras hold, i.e.
\begin{align*}
&[\so(V)_\pm,\so(V)_\pm]\subset \so(V)_+,\qmq{and}[\so(V)_+,\so(V)_-]\subset \so(V)_-\;.
\end{align*}
For an isometric immersion $f:M\to N$ and some $p\in M$ we will apply this
construction with
$W=T_pM$ and $V=\osc_pf$ resp.\ $V=T_{f(p)}N$\,; then we obtain the induced splitting
\begin{equation}\label{eq:so_osc_2}
\so(T_{f(p)}N)=\so(T_{f(p)}N)_+\oplus\so(T_{f(p)}N)_-\qmq{and}\so(\osc_pf)=\so(\osc_pf)_+\oplus\so(\osc_pf)_-\;.
\end{equation}

\bigskip
\begin{definition}\label{de:fetth}
For each $p\in M$ let $\fetth:T_pM\to\so(T_{f(p)}N)$ be the linear map
defined by 
\begin{equation}\label{eq:fetth}
\forall\,x,y\in T_pM\,,\xi\in\bot_pM:\fetth(x)(y+\xi):=h(x,y)-S_\xi x\;;
\end{equation}
note that $h$ and $\fetth$ are equivalent objects.
\end{definition}

In the following,
$\so(\osc_pf)$ is seen as a Lie subalgebra of $\so(T_{f(p)}N)$ in a natural way:
\begin{equation}\label{eq:so_osc}
\so(\osc_pf)\cong\Menge{A\in \so(T_{f(p)}N)}{A(\osc_pf)\subset\osc_pf,\;A|(\osc_pf)^\bot=0}\;;
\end{equation}
then we have 
\begin{align}\label{eq:split_osc}
&\so(\osc_pf)_\pm=\so(T_{f(p)}N)_\pm\cap\so(\osc_pf)\;,\\
\label{eq:fetth_in_osc}
&\forall x\in T_pM:\;\fetth(x)\in\so(\osc_pf)_-\;.
\end{align} 
For a submanifold $M\subset N$ the linear spaces
$\so(\osc_pM)_\pm\subset\so(T_pN)_\pm$ and the linear map $\fetth:T_pM\to\so(\osc_p M)_-$ are
defined by means of $\iota^M$\,. 

Throughout this paper, we will make use of the following convention: Given two
points $p,q\in M$ and a linear map $\ell:T_{f(p)}N\to T_{f(q)}N$ with
$\ell(\osc_pf)\subset \osc_qf$ we put 
\begin{equation}\label{eq:A_osc}
\ell^{\osc}:=\ell|\osc_pf:\osc_pf\to\osc_qf\;.
\end{equation} 

\bigskip
\begin{theorem}\label{th:hol}
Let $f:M\to N$ be a parallel isometric immersion defined on a symmetric space
$M$\,. The extrinsic holonomy Lie algebra of $\osc f$ is characterized
by the following properties:
\begin{enumerate}
\item There is the splitting
\begin{equation}\label{eq:splitting_of_hol}
\hol(\osc f)=\hol(\osc f)_+\oplus \hol(\osc f)_-\;,
\end{equation}
with $\hol(\osc f)_\pm:=\hol(\osc f)\cap\so(\osc_of)_\pm$\,.
\item 
We have $R^N(x,y)\,\osc_of\subset\osc_of$ and
$R^N(\xi,\eta)\,\osc_of\subset\osc_of$ for all $x,y\in T_oM\,,\xi,\eta\in\bot^1_of$\,,
and the splitting~\eqref{eq:splitting_of_hol} is given by
\begin{align}\label{eq:hol_plus}
&\hol(\osc f)_+=\Spann{\big(R^N(x,y)\big)^{\osc}}{x,y\in T_oM}+\Spann{\big (R^N(\xi_1,\xi_2)\big )^{\osc}}{\xi_1,\xi_2\in\bot^1_of}\;,\\
\label{eq:hol_minus}
&\hol(\osc f)_-=\Spann{[\fetth(x),A]}{x\in T_oM, A\in\hol(\osc f)_+}\;.
\end{align}
\end{enumerate}
Furthermore, for all $x\in T_oM$ we have
\begin{equation}\label{eq:[fetth(x),_]}
[\fetth(x),\hol(\osc f)]\subset\hol(\osc f)\;.
\end{equation}
If moreover $f$ is a full immersion, then $\hol(\osc f)\cong\hol(f^*TN)$\,,
more precisely:
\begin{enumerate}
\addtocounter{enumi}{2}
\item 
We have $A(\osc_of)\subset\osc_of$ for all $A\in\hol(f^*TN)$\,, and the linear map
 $
\hol(f^*TN)\to\hol(\osc f),\;A\mapsto A^{\osc f}
$
is a Lie algebra isomorphism.
\end{enumerate}
\end{theorem}
The proof of Theorem~\ref{th:hol} is given in Section~\ref{s:hol}.

\bigskip
In Section~\ref{s:extrinsically_hol}, for every full symmetric submanifold $M$
of some simply connected symmetric spaces $N$ the extrinsic holonomy Lie algebra 
is calculated in an explicit way. Because of the following result, thereby it is always enough to consider the
case when $N$ is an irreducible Riemannian space:

\bigskip
\begin{theorem}\label{th:Naitoh_2}{\cite{N3}}
Let $N$ be a simply connected symmetric space, $N\cong R^d\times
N_1\times\cdots\times N_k$ its ``de\,Rham decomposition'' (see~\cite{KN},~Ch.\,IV,
Theorem~6.2) and $M\subset N$ a
symmetric submanifold. 
Then there exist symmetric submanifolds $M_0\subset R^d$ and $M_i\subset N_i$
for $i\geq 1$ such that $M\cong M_0\times M_1\times\cdots\times M_k$ (as a submanifold). 
\end{theorem} 
In the irreducible case, one knows  the following result, which is a consequence of~Proposition 9.3.3 combined with Theorem
9.3.4 from~\cite{BCO}:

\bigskip
\begin{theorem}[Naitoh]\label{th:Naitoh}
\item If $N$ is a simply connected, irreducible symmetric space and
$M\subset N$ is a full symmetric submanifold with $o\in M$\,, then only the
following possibilities can occur:\footnote{As was shown
  in~\cite{Ko},~\cite{NT} and~\cite{Ts4}, actually there do not exist any full parallel submanifolds in a complex or quaternionic hyperbolic space.}
\begin{itemize}
\item $N$ is a real space form.
\item $N^{2n}$ is a complex space form ($n\geq 2$) and $M$ is a complex submanifold.
\item $N^{2n}$ is a complex space form ($n\geq 2$) and $M^n$ is a Lagrangian submanifold.
\item $N^{4n}$ is a quaternionic space form ($n\geq 2$) and $M^{2n}$ is a totally complex
  submanifold.
\item  The rank of $N$ is larger than 1, $N$ admits a symmetric R-space (see
  Definition~\ref{de:symmetric_R-space}) and $M$ belongs to the family of symmetric
  submanifolds associated therewith (in the sense of Definition~\ref{de:BENT}). 
\end{itemize}
\end{theorem}

We will prove:
 
\bigskip
\begin{theorem}\label{th:hol_for_symmetric_M}
Let $N$ be a simply connected, irreducible symmetric space and $M^m$ be a
full symmetric submanifold of $N$ through $o$ with $m\geq2$\,. The extrinsic holonomy Lie
algebra of $M$ is given as follows:
\begin{enumerate}
\item Suppose that $N$ is a Hermitian symmetric space. Here we have
  $\hol(N)=[\hol(N),\hol(N)]\oplus\R\,\j$\,, where $\j$ denotes the complex
  structure of $N$ at $o$ and $[\hol(N),\hol(N)]$ is the commutator ideal of
$\hol(N)$\,. If $M$ is a Lagrangian submanifold of $N$\,, then we have $\hol(TN|M)=[\hol(N),\hol(N)]$\,. 
\item Suppose that $N$ is the quaternionic projective space $\bbH\rmP^n$
with $n\geq 2$\,. Here we have $\hol(N)=\sp(T_oN)\oplus\scrQ$\,,
where $\scrQ\subset\so(T_oN)$ denotes the quaternionic structure at $o$\,. 
For each complex submanifold $M^{2n}\subset N$ (in the sense
of~\cite{Ts4}) it is possible to choose a canonical basis $\{\i,\j,\k\}$
of $\scrQ$ (i.e. $\{\i,\j,\k\}$ is a basis of $\scrQ$ such that the usual quaternionic relations 
$\i^2=\j^2=-\Id$\,, $\i\circ\j=-\j\circ\i=\k$ hold, see~\cite{Ts4}, Definition~2.2)
such that $\i(T_oM)=T_oM$ and $\j(T_oM)=\bot_oM$ holds. In this situation, we have $\hol(TN|M)=\sp(T_oN)\oplus\R\,\i$\,. 
\item Suppose that $M^2$ is a submanifold of a 4-dimensional space form
(this ``exceptional'' case actually occurs, see Remark~\ref{re:exception} below). Then we have
$\hol(N)\cong\so(4)$\, and either $\hol(TN|M)=\hol(N)$ or
$\dim_\R(\hol(TN|M))=4$ holds\,.\footnote{In this case, the result is not very satisfying. However, the methods developed in this paper
are not suitable to obtain a better result. Possibly, here a ``case by case'' argument would shed some more light on the subject.} 
\item In all other cases we have $\hol(TN|M)=\hol(N)$\,.
\end{enumerate}
Moreover, if in Case~(a) the rank of $N$ is
larger than 1, then there exists some $x\in T_oM$ with
$\fetth(x)=\j$\,.\footnote{Hence we see that in this case the condition $\fetth(T_oM)\subset\hol(TN|M)$ is violated.} 
\end{theorem}
The proof of Theorem~\ref{th:hol_for_symmetric_M} can be found in Section~\ref{s:extrinsically_hol}.
Let us ``apply'' Theorem~\ref{th:hol_for_symmetric_M} to the relevant cases
(actually, Theorem~\ref{th:hol_for_symmetric_M} gets proved the other way around):

\bigskip
\begin{example}
In the following we assume that $n,m\geq 2$\,.
\begin{enumerate}
\item
Full symmetric submanifolds $M^m$ of the Euclidian sphere $S^n$  resp.\ of the hyperbolic space $\rmH^n$ where classified
in~\cite{Fe0} resp.\ in~\cite{BR} and~\cite{Ta}. In accordance with
Theorem~\ref{th:hol_for_symmetric_M}\,, here we have
$\hol(TN|M)=\hol(N)=\so(T_oN)$ unless $(n,m)=(4,2)$\,.
\item 
Full, complex symmetric submanifolds $M\subset\C\rmP^n$ where classified
in~\cite{NaTa}; cf.\ also Table~9.1 of~\cite{BCO}. A prominent example is given by the image of the ``Segre embedding'' $f:\C\rmP^m\times
\C\rmP^n\to\C\rmP^N$ with $N+1=(m+1)(n+1)$\,, given by
$$
([z_0:\cdots :z_m],[w_0:\cdots :w_n])\mapsto [z_0w_0:z_0w_1:\cdots
:z_mw_n]\ \text{(all possible combinations)}\;.
$$
In accordance with Theorem~\ref{th:hol_for_symmetric_M}\,, here we always have $\hol(TN|M)=\hol(N)=\fraku(T_oN)$\,. 
\item 
Full, Lagrangian symmetric submanifolds $M\subset \C\rmP^n$ where classified
in~\cite{N}; cf.\ also Table 9.2
of~\cite{BCO}. In accordance with Theorem~\ref{th:hol_for_symmetric_M}, here
we always have $\hol(TN|M)=\su(T_oN)$\,, which is strictly contained in $\hol(N)$\,.
\item
Full, totally complex symmetric submanifolds $M^{2n}\subset\bbH\rmP^n$ where classified
in~\cite{Ts4}; cf.\ also Table 9.4 of~\cite{BCO}.
\item
Let $N$ be an irreducible Hermitian symmetric space % Then we have the splitting $\frakk=\frakc\oplus[\frakk,\frakk]$ with a
% one-dimensional factor $\frakc$\,, the center of $\frakk$
% (cf.~\cite{BCO},~A.4); furthermore, there exists a canonical isomorphism $\hol(N)\cong\frakk$ (see Proposition~\ref{p:frakk_gleich_hol}). 
which admits a symmetric R-space and let $\{M_c\}$ (where
$c$ ranges over $\R$) denote the family of symmetric submanifolds associated
therewith, see Definition~\ref{de:BENT}. As a consequence of Proposition~\ref{p:BENT}, $M_c$ is a full submanifold of $N$ unless
$c=0$\,. Moreover, $M_c$ is a Lagrangian submanifold
of $N$ for each $c\in \R$\,, by virtue of Lemma~\ref{le:frakc}; therefore, in
accordance with Theorem~\ref{th:hol_for_symmetric_M}, the extrinsic holonomy
Lie algebra $\hol(TN|M_c)$ is given by $[\hol(N),\hol(N)]$ for each $c\neq 0$\,. 
\item
Let $N$ be an irreducible symmetric space which is not of Hermitian
type and admits a symmetric
R-space. As above, let $\{M_c\}$ denote the family of
symmetric submanifolds associated therewith; then  $\hol(TN|M_c)=\hol(N)$ for each $c\neq 0$\,.
\end{enumerate}
\end{example}

\bigskip
\begin{remark}\label{re:exception}
In $S^4$ there exists a symmetric submanifold $M^2$ which is isometric to $\R\rmP^2$\,;
it is congruent to a standard embedded symmetric R-space which is a (symmetric) orbit 
of the isotropy representation of the five dimensional symmetric space
$\SU(3)/SO(3)$\,. By the results of~\cite{BR} and~\cite{Ta}, one knows that there also exist
certain full symmetric submanifolds $M^2$ in $\rmH^4$\,. 
More precisely, $M$ is an extrinsic product (in the sense
of~\cite{BR},~Definition~4) 
which is isometric to $\R\times S^1(r)$\,, $S^1(r)\times S^1(s)$ or $S^2(r)$ 
(where $r,s>0$ are arbitrary). 
\end{remark}

As a conclusion of Theorem~\ref{th:hol_for_symmetric_M}, we notice that \textbf{for every full
symmetric submanifold $\mathbf{M}$ of some simply connected, irreducible symmetric space
$\mathbf{N}$ the subspace $\mathbf{\hol(TN|M)\subset \hol(N)}$ is
surprisingly large; always its codimension is 0, 1 or 2.} 
Moreover, Parts~(c) and~(d) of
Lemma~\ref{le:about_hol} in Section~\ref{s:extrinsically_hol} show that in
each case $\hol(TN|M)$ is ``as large as possible'' 
(maybe the last assertion is not true for the exceptional case described in Part~(c) of
Theorem~\ref{th:hol_for_symmetric_M}).
\section{Symmetric submanifolds}
\label{s:results}
We aim to review the relation between parallel
and symmetric submanifolds of a symmetric space $N$\,.
Let $\Iso(N)$ denote the Lie group of isometries on $N$ (see~\cite{He},~Ch.\,IV,~\S~2\ \text{and}~\S~3), $\Iso^0(N)$ its connected
component and $\fraki(N)$ the corresponding
Lie algebra. For each $X\in\fraki(N)$ we have the one-parameter subgroup
$\psi_t^X:=\exp(t\,X)$ of isometries on $N$\,; the corresponding ``fundamental vector
field'' $X^*$ on $N$ (in the sense of~\cite{KN}) defined by
\begin{equation}\label{eq:fundamental_vector_field}
X^*(p):=\frac{\diff}{\diff t}\Big|_{t=0}\psi_t^X(p)
\end{equation}
is a Killing vector field on $N$ such that $\psi_t^X$ ($t\in\R$) is the flow
of $X^*$\,.\footnote{By the map $X\mapsto X^*$ the vector space $\fraki(N)$ is identified with the Lie
algebra of Killing vector fields on $N$\,; but one should be aware that
$[X^*,Y^*]=-[X,Y]^*$\,, where the bracket on the l.h.s.\ is the Lie bracket for
vector fields and on the r.h.s.\ is the bracket of $\fraki(N)$\,.}
The {\em isotropy subgroup}\ of $\Iso^0(N)$ at some fixed origin $o\in $ is by definition 
\begin{equation}\label{eq:rmK}
\rmK:=\Menge{g\in\Iso^0(N)}{g(o)=o}\;.
\end{equation}

The {\em isotropy
representation}\ is given by 
$$
\rmK\to\SO(T_oN), g\mapsto T_og\;.
$$
Let $\frakk$ denote the Lie algebra of $\rmK$ and $\pi_2:\frakk\to\so(T_oN)$ the linearized isotropy representation,
i.e.
\begin{equation}\label{eq:Kovv_X_1}
\forall X\in\frakk,u\in T_oN:\pi_2(X)=\frac{\diff}{\diff t}\Big|_{t=0}T_o\psi_t^X(u)\;.
\end{equation}

\bigskip
\begin{theorem}[Str\"ubing-Naitoh-Eschenburg]\label{th:struebing}
For a submanifold $M$ of a simply connected symmetric space $N$ the following
assertions are equivalent:
% \Abstand{3}
\begin{enumerate}
\item
$M$ is a symmetric submanifold\,.
\item
$M$ is a complete parallel submanifold\,, such that all
normal spaces are curvature invariant.
\item $M$ is a complete, parallel submanifold,, at one point $p\in M$
  the normal space $\bot_pM$ is curvature invariant, and
 \begin{equation}\label{eq:Eschenburg}
\fetth(T_pM)\subset\pi_2(\frakk)\;.
\end{equation}

\end{enumerate}
\end{theorem}
For a proof of $(a)\Leftrightarrow(b)$ see~\cite{N2}, Corollary~1.4, for the
other directions see~\cite{E2},~Theorem~4.
Notice that if $N$ has constant
curvature, then every subspace of $T_pN$ is curvature invariant;
therefore, as a consequence of Theorem~\ref{th:struebing}, in a space form
{\em every}\/ complete parallel submanifold is extrinsically symmetric; and
the converse is also true:

\bigskip
\begin{proposition}\label{ex:constant_curvature}
If every complete parallel submanifold of $N$ is extrinsically symmetric, then $N$ is of constant curvature.
\end{proposition}
\begin{proof}
A particular example for parallel submanifolds are geodesic lines, see
Example~\ref{ex:circles}.
Suppose that all complete geodesic lines of $N$ would be
extrinsically symmetric. Then according to Theorem~\ref{th:struebing}
for each $p\in M$ any linear hyperplane of $T_pN$ is curvature invariant
(since it can be realized as the normal space of some geodesic through $p$).
By a result due to E.~Cartan, $N$ is a space of constant curvature, see~\cite{Ts2}.
\end{proof}
%  a higher dimensional, complete, parallel
% submanifold which is neither extrinsically symmetric nor full, and there
% also exists a full, extrinsically homogenous (i.e. obtained as an orbit of a
% subgroup of the isometry
% group of the ambient space), complete parallel submanifold which is not extrinsically symmetric:

\bigskip
\begin{example}
\begin{enumerate}
\item Let $\R\rmP^n$ be canonically embedded in $\C\rmP^n$ and let $M$ be a
(proper) symmetric submanifold of $\R\rmP^n$\,. Then $M$
is parallel in $\C\rmP^n$\,, but not extrinsically symmetric in $\C\rmP^n$\,.
\item The circle mentioned in Example~\ref{ex:CP2} is a covering
  onto a full parallel submanifold $M\subset \C\rmP^2$\,. Moreover, $M$ is the orbit of a subgroup
  of $\Iso(\C\rmP^2)$\,, but $M$ is not a symmetric submanifold of $\C\rmP^2$\,.
\end{enumerate}
\end{example}
\begin{proof}
For~(a): Of course, $M$ is also parallel in $\C\rmP^n$\,.
On the other hand, the normal spaces $\bot_pM\;(p\in M)$ are
neither complex nor totally real subspaces. It is well known that therefore they are not
curvature invariant. Thus $M$ is not extrinsically symmetric in $\C\rmP^n$ as
a consequence of Theorem~\ref{th:struebing}.

For~(b): By a result of~\cite{MT}, $c$ is the orbit of a one-parameter
subgroup of $\Iso(\C\rmP^2)$ and hence $c$ is a covering onto a
full, extrinsically homogenous parallel submanifold of $\C\rmP^2$\,, which can be
not be extrinsically symmetric in $\C\rmP^2$ according to Example~\ref{ex:CP2} in
combination with Theorem~\ref{th:struebing}.
\end{proof}

\subsection{Irreducible symmetric R-spaces}
\label{s:symmetric_R-space} 
For this section cf.~\cite{BCO}, Ch.~3.7 and~A.4\,,~\cite{BENT} and~\cite{EH}.
Let $N$ be a simply connected, irreducible symmetric space; hence $N$ is of
compact type or of non-compact type. Let $o\in N$ be some origin, $\rmK\subset\Iso^0(N)$ the isotropy
subgroup, $\fraki(N)=\frakk\oplus\frakp$ the Cartan decomposition and $B$
the Killing form of $\fraki(N)$\,. We consider the adjoint
representation $\Ad:\Iso(N)\to\Gl(\fraki(N))$ and its linearization
$\ad:\fraki(N)\to\gl(\fraki(N))$\,. Let $\epsilon\in\{-1,1\}$ be chosen such
that the restriction of $\epsilon\,B$ to $\frakp\times\frakp$ is a positive
definite inner product; hence $\epsilon=1$ if and only if $N$ is of
non-compact type. Then $\Ad$ induces a faithful orthogonal representation of
$K$ on $\frakp$\,, by restriction; the corresponding infinitesimal action is
given by $\ad_\frakp:\frakk\to\so(\frakp)$\,. In this section, we will consider certain $\Ad(\rmK)$-orbits of $\frakp$\,, 
so called standard embedded irreducible {\em symmetric R-spaces}\/. As was shown
in~\cite{Fe0}, these objects are the fundamental examples of parallel
submanifolds in a Euclidian space; moreover, they also give rise to families
of symmetric submanifolds in $N$ as will be explained in the next section.

\bigskip
\begin{remark}\label{re:symmetric_R-space}
Let $N^*$ denote the dual symmetric space (cf.~\cite{BCO},~A.4), which again is a simply connected
irreducible symmetric space such that $\epsilon^*=-\epsilon$\,. Then $K$ is also the isotropy group of $N^*$ 
and $\fraki(N^*)=\frakk^*\oplus\frakp^*\cong\frakk\oplus\i\,\frakp$ (seen as a Lie subalgebra of
the complexified Lie algebra $\fraki(N)\otimes\C$) is the Cartan decomposition for $N^*$\,. 
\end{remark} 

\bigskip
\begin{lemma}\label{le:ad(X)}
For each $X\in\frakp$ with $\ad(X)^3=\epsilon\,\ad(X)$ 
%put $\sigma:=\exp(\pi\,X)\in\Iso^0(N)$\,. Then 
we have:
\begin{enumerate}
\item $\ad(\i\,X)^3=-\epsilon\,\ad(\i\,X)$ on $\frakk\oplus\i\,\frakp$ (in the sense of Remark~\ref{re:symmetric_R-space}).
\item $\ad(X)$ is diagonalizable over $\C$ with $\Spec(\ad(X))\subset\{-1,0,1\}$
    (in case $\epsilon=1$) resp. $\Spec(\ad(X))\subset\{-\i,0,\i\}$ (otherwise). 
Hence, there exist decompositions
\begin{align}
\label{eq:decomposition_of_frakk_with_respect_to_ad(X)}
&\frakk=\frakk_0\oplus\frakk_\epsilon:=\Kern(\ad(X)|\frakk)\oplus\Menge{Y\in\frakk}{\ad(X)^2Y=\epsilon\,Y}\;,\\
\label{eq:decomposition_of_frakp_with_respect_to_ad(X)}
&\frakp=\frakp_0\oplus\frakp_\epsilon:=\Kern(\ad(X)|\frakp)\oplus\Menge{Y\in\frakp}{\ad(X)^2Y=\epsilon\,Y}\;.
\end{align}
Then~\eqref{eq:decomposition_of_frakp_with_respect_to_ad(X)} is an orthogonal
splitting and we have
\begin{equation}\label{eq:ad_frakp}
\ad_\frakp(\frakk_0)\subset\so(\frakp)_+\qmq{and}\ad_\frakp(\frakk_\epsilon)\subset\so(\frakp)_-\,.
\end{equation} 

\item We have $\ad(X)\,\frakk_\epsilon\subset\frakp_\epsilon$\,,
  $\ad(X)\,\frakp_\epsilon\subset\frakk_\epsilon$ and $J:=-\epsilon\,\ad(X)|\frakk_\epsilon\oplus\frakp_\epsilon$\,,
  seen as an endomorphism of $\frakk_\epsilon\oplus\frakp_\epsilon$\,,
  satisfies $J^2=\epsilon\,\id$\,; moreover,
\begin{equation}
\label{eq:ad(X)_2}
J|\frakp_\epsilon:\frakp_\epsilon\to\frakk_\epsilon\qmq{\ is a linear isomorphism\;.}
\end{equation}
\end{enumerate}
\end{lemma}
\begin{proof}
(a) is obvious. Therefore, by switching between $N$ and the dual symmetric space $N^*$\,, we may assume that $\epsilon=1$ (i.e. $N$ is of non-compact
type); then the other results are in accordance with~\cite{BCO},~Example~7.7.
\end{proof}
% We will see $\fraki(N)$ as a Euclidian vector space equipped with the inner
% inner product $B_\sigma(X,Y):=-B(X,\sigma Y)$\,, where $\sigma$ denotes the Cart. Thus $\ad(X)$ is a  skew-symmetric endomorphism on $\fraki(N)$\,, hence diagonalizable over $\C$\,.
% Moreover, the complex spectrum of $\ad(X)$ is contained in $\{0,\i,-\i\}$\,, with $\i:=\sqrt{-1}$\,.
% Using also the bracket relations for the Cartan decomposition,
% Equations~\eqref{eq:decomposition_of_frakk_with_respect_to_ad(X)}-\eqref{eq:ad(X)_2}
% follow easily. 
% Hence we have for all $t\in\R$\,:
% \begin{align}\label{eq:Ad(exp(tX))_1}
% &\forall Y\in\frakk_+\oplus\frakp_+:\exp(t\,\ad(X)\big)\,Y=Y\qmq{and}\\
% \label{eq:Ad(exp(tX))_2}
% &\forall Y\in\frakk_-\oplus\frakp_-:\exp(t\,\ad(X))\,Y=\cos(t)\,Y+\sin(t)J\,Y\;.
% \end{align}
% Thus $\frakk_+\oplus\frakp_+$ and $\frakk_-\oplus\frakp_-$ are the
% $+\,1$-eigenspace and the $-1$-eigenspace of
% $\Ad(\sigma)$\,, respectively. Equation~\eqref{eq:ad_frakp} now follows, because $\sigma$
% respects the brackets.

% Because of Remark~\ref{re:symmetric_R-space}, our definition of symmetric
% R-spaces stands in accordance with~\cite{BCO},~A.\,4. 

\bigskip
\begin{proposition}\label{p:symmetric_R-space}
Let $N$ be a simply connected symmetric space of non-compact type, suppose
that there exists $X\in\frakp$ with $X\neq 0$ and $\ad(X)^3=\epsilon\,\ad(X)$ and consider the orbit
$M:=\Ad(\rmK)\,X$\,. Then we have  $T_XM=\frakp_\epsilon$\,,
$\bot_XM=\frakp_0$ and
\begin{equation}\label{eq:fetth_1}
\forall Y\in\frakp_\epsilon:\fetth(Y)=\ad(J\,Y)|\frakp:\frakp\to\frakp\;,
\end{equation}
where $J$ is the linear map from Lemma~\ref{le:ad(X)}. Moreover, $M$ is a 1-full symmetric
submanifold of $\frakp$ and $h_o$ is non-degenerate\,.

\end{proposition}
\begin{proof} 
% ``extrinsic symmetry'': Because $M$ is an orbit under some isometric action on $\frakp$\,, it will follow
% that $M$ is a symmetric submanifold of $\frakp$\,, if we can show that
% $\Ad(\sigma)|\frakp:\frakp\to\frakp$ is the extrinsic symmetry at the point $X\in M$
% in accordance with Definition~\ref{de:extrinsically_symmetric}.
We have
$$T_XM=\ad(\frakk)\,X=[\frakk,X]\stackrel{\eqref{eq:decomposition_of_frakk_with_respect_to_ad(X)}}=[\frakk_\epsilon,X]
\stackrel{\eqref{eq:ad(X)_2}}=\frakp_\epsilon\;.$$
Thus $T_XM=\frakp_\epsilon$ and hence $\bot_XM=\frakp_0$\,, since the
splitting~\eqref{eq:decomposition_of_frakp_with_respect_to_ad(X)} is an
orthogonal sum.

For~\eqref{eq:fetth_1}:  
For each $Y\in\frakp_\epsilon$ we have $J\,Y\in\frakk_\epsilon$ according
to~\eqref{eq:ad(X)_2}; hence the linear map $A_Y:=\ad_\frakp(JY)$ (seen as
a linear vector field on $\frakp$) is tangent to $M$\,.
We have (by the Gau{\ss} equation) for all $Z\in \frakp_\epsilon$
$$
\fetth(Y)\,Z=\epsilon\,h(\,J^2Y,Z)=\epsilon\,h([X,[X,Y]],Z)=-\epsilon\,h(A_YX,Z)=-\epsilon\,(A_YZ)^\bot
=[J\,Y,Z]^\bot\stackrel{~\eqref{eq:ad_frakp}}=[J\,Y,Z]\;;
$$
therefore~\eqref{eq:fetth_1} follows in view of~\eqref{eq:odd},~\eqref{eq:ad_frakp}.
% For this: We have $\Ad(\sigma)\,X=X=\Ad\exp(\pi\,X)\,X=X$\,,
% and since furthermore~\eqref{eq:decomposition_of_frakp_with_respect_to_ad(X)} gives the
% decomposition $T_XN=T_XM\oplus\bot_XM$\,, we deduce from
% Lemma~\ref{le:ad(X)} that $\Ad(\sigma)|\frakp:\frakp\to\frakp$ is the reflection
% in the affine subspace $\bot_XM$ of $\frakp$\,.
% On the other hand, $\Ad(\sigma)\frakk=\frakk$ and hence $\sigma K\sigma=K$\,, because $K$
% is connected. Thus we have for each $k\in K$\,:
% $$
% \Ad(\sigma)(\Ad(k)\,X)=\Ad(\sigma k)\,X=\Ad(\sigma
% k\sigma)\Ad(\sigma)\,X=\Ad(\sigma k\sigma)\,X\in \Ad(\rmK)\,X\;,
% $$
% and hence $\Ad(\sigma)(M)=M$\,.
Now the non-degeneracy of $h$ follows, since $J|\frakp_\epsilon$ in injective
and $\ad_\frakp$ is a faithful representation. Furthermore, it is well known that $M$ is a symmetric submanifold of
$\frakp$\,, cf~\cite{BCO},~Proposition~3.7.7 or~\cite{EH}, Theorem~2.
To see that $M$ is 1-full in $\frakp$\,, let $Z\in\bot_XM=\frakp_0$ be given. I claim
that $S_Z=0$ already implies that $Z=0$\,:

Thereby, without loss of generality we may assume that $\epsilon=1$\,. According
to~\eqref{eq:decomposition_of_frakp_with_respect_to_ad(X)}, we have
$[X,Z]=0$\,. Therefore we may choose a maximal Abelian subspace $\fraka\subset\frakp$ with $\{X,Z\}\subset
\fraka$\,. Then the adjoint action of $\fraka$ on $\fraki(N)$ is simultaneously diagonalizable. 
Let $\Sigma$ denote the corresponding set of weights, usually called the
``restricted roots'', choose some ordering of $\fraka$ such that $X$ lies in
the closure of the Weyl Chamber where the roots are positive 
(cf.~\cite{He},~Ch.\,VII, Remark after Lemma~2.20) 
and let $\Sigma^+$ denote the set of positive roots. Put
$\Sigma^+_n:=\Menge{\lambda\in\Sigma^+}{\lambda(X)=n}$ for $n=0,1$\,;
then, since $\Spec(\ad(X))\subset\{-1,0,1\}$\,, $\Sigma^+$ is
the disjoint union of $\Sigma^+_0$ and $\Sigma^+_1$\,.
In compliance with a result of~\cite{BCO} (see p.~64 there), the set of eigenvalues for $S_Z$ is given by
$\Menge{\lambda(Z)}{\lambda\in \Sigma^+_1}$\,.
Therefore, if $S_Z=0$\,, then $\lambda(Z)=0$ for all $\lambda\in\Sigma^+_1$\,. 
I claim that this already implies that $\lambda(Z)=0$ for all
$\lambda\in\Delta^+$\,: 

For this, note that $\Sigma^+_0$ is the intersection
of $\Sigma^+$ with  the hyperplane
$\Menge{\lambda\in\fraka^*}{\lambda(X)=0}$\,. Therefore, since $\Sigma$ is a
root system, we see that $\Sigma^+_1$ spans a
vector subspace $V\subset \fraka^*$ which is invariant under all reflections in the
various elements of $\Sigma^+$\,. Since the (abstract) Weyl group of $\Sigma$ is generated by the
reflections in the various elements of $\Sigma^+$\,, we hence conclude from
the irreducibility of $N$ that $V=\fraka^*$ holds, which immediately gives
our claim.

Hence $\lambda(Z)=0$ for all
$\lambda\in\Delta^+$\,. To see that this already implies the vanishing of $Z$\,, we proceed as
follows: Let $\frakh$ be a maximal Abelian
subspace of $\frakg^*$ with $\fraka\subset\frakh$ and let $\frakh_\C$ denote its complexification;
then $\frakh_\C$ is a Cartan subalgebra of $\frakg_\C$\,, according to~\cite{He},~Ch.\,VI,~Lemma~3.2. Let
$\Delta$ denote the corresponding set of roots; hence $\lambda:\frakh_\C\to \C$ is a linear
function with $\lambda(\frakh)\subset\R$ for each $\lambda\in\Delta$\,. As
explained in~\cite{He},~Ch.\,VI,~\S~3, for each $\lambda\in\Delta$ we either
have $\lambda|\fraka=0$\,, or $\lambda|\fraka\in\Sigma$ holds; therefore, since $\Delta$
spans the dual space $\frakh_\C^*$ (cf.~\cite{He},~Ch.\,III,~\S~4, Proof of
Theorem~4.2, Equation~(2)), we obtain that $\lambda(Z)=0$ for all
$\lambda\in\Sigma$ implies $Z=0$\,. 

Thus $\bot_XM\to\End(T_XM)\,,Z\mapsto S_Z$ is an injective map, and hence a
straight forward calculation shows
that the second fundamental form of $M$ spans $\bot_XM$\,. This finishes the proof.
\end{proof}

\bigskip
\begin{definition}\label{de:symmetric_R-space}
Let $N$ be a simply connected, irreducible symmetric space of compact type or of non-compact type, $o\in N$ some origin, $K\subset\Iso^0(N)$ the
isotropy group and $\fraki(N)=\frakk\oplus\frakp$ the Cartan decomposition.
If there exists $X\in\frakp$ with $X\neq 0$ and $\ad(X)^3=\epsilon\,\ad(X)$\,, 
then $M:=\Ad(\rmK)\,X$ is called a (standard embedded, irreducible) {\em
  symmetric R-space}\/ and we say that ``$N$ admits a symmetric R-space''.
\end{definition}
%Irreducible symmetric R-space are sometimes also called ``(irreducible) symmetric real flag manifolds''.

Because of Remark~\ref{re:symmetric_R-space} and Lemma~\ref{le:ad(X)}, on the level of symmetric
R-spaces, it is always enough to consider the case when $N$ is of compact
type, i.e. $\epsilon=-1$\,.
Then the following theorem gives the classification of irreducible symmetric
R-spaces, in accordance with~\cite{BCO},~Tables~A.6 and A.7\,:

\bigskip
\begin{theorem}[\cite{KNa}]
\label{th:The_classification_of_irreducible_symmetric_R-spaces}
The irreducible symmetric spaces $N$ which are of compact type and admit a
symmetric R-space $M$ are given as follows:
\begin{enumerate}
\item
If additionally $N$ is of Hermitian type

\begin{tabular}{c|c|c}\label{t:Hermitian}
$N$                      & $M$              & Remarks\\ \hline    
$\SU(2n)/\rmS\big(\rmU(n)\times\rmU(n)\big)$ & $\rmU(n)$                  &  $n\geq 2$\\
$\SO(n+2)/\SO(2)\times\SO(n)$                & $\big(\rmS^1\times\rmS^{n-1}\big)/\Z_2$      &  $n\geq 3$\\
$\Sp(n)/\rmU(n)$                                & $\rmU(n)/\SO(n)$    & $n\geq 3$\\
$\SO(4n)/\rmU(2n)$                           & $\rmU(2n)/\Sp(n)$  & $n\geq 3$\\ 
$E_7/\rmT\cdot E_6$                             & $(\rmT\cdot E_6)/F_4$      &$--$\\
\end{tabular}

\item
Otherwise

\begin{tabular}{c|c|c}
 $N$                                         & $M$        &  Remarks\\ \hline              
 $\mathrm{Spin}(n)$                          & $\SO(n)/\big(\SO(2)\times\SO(n-2)\big)$         & $n\geq 5$\\
 $\mathrm{Spin}(2n)$                         & $\SO(2n)/\rmU(n)$       & $n\geq 3$\\
 $\SU(n)$                                    & $\SU(n)\big/\rmS(\rmU(p)\times\rmU(n-p))$  & $n\geq 2$\,, $1\leq p\leq [\frac{n}{2}]$\\
 $\Sp(n)$                                    & $\Sp(n)/\rmU(n)$                       & $n\geq 2$\\
 $\rmE_6$                                    & $\rmE_6\big/\rmT\cdot\mathrm{Spin}(10)$  & $--$\\
 $\rmE_7$                                    & $\rmE_7\big/\rmT\cdot \rmE_6$          & $--$\\
$\SU(n)/\SO(n)$                              & $\rmG_p(\R^n)$
 & $n\geq 3$\,, $1\leq p\leq[\frac{n}{2}]$\\
$\SU(2n)/\Sp(n)$                             & $\rmG_p(\bbH^n)$       & $n\geq 2$\,, $1\leq p\leq[\frac{n}{2}]$\\
$\SO(2n)\big/\SO(n)\times \SO(n)$            & $SO(n)$                 & $n\geq 5$\\
$\Sp(2n)\big/\Sp(n)\times \Sp(n)$            & $Sp(n)$                & $n\geq 2$\\
$\rmE_6/\Sp(4)$                              & $\rmG_2(\bbH^4)/\Z_2$  & $--$\\
$\rmE_6/\rmF_4$                              & $\bbO P^2$            & $--$\\
$\rmE_7/\SU(8)$                              &$\big(\SU(8)/\Sp(4)\big)\big/\Z_2$ & $--$\\
$\SO(n)\big/\SO(p)\times\SO(n-p)$
&$(\rmS^{p-1}\times\rmS^{n-p-1})/\Z_2$  & $n\geq 3$\,, $3\leq p\leq
[\frac{n}{2}]$\\ 
\end{tabular}
\end{enumerate}
\end{theorem}
% \bigskip
% The following result on the classification of parallel submanifolds in a euclidian space
% is is in accordance with~Theorem~3.7.8 in~\cite{BCO}, see also Theorem~2 in~\cite{EH}:
% \bigskip
% \begin{theorem}\label{th:Ferus}[Ferus]
% \label{th:extrinsically_symmetric}
% For a full, parallel
% submanifold $M\subset \R^n$ there exists a symmetric space $N$ of
% compact type,
% \item 
% $X\in\frakp$ with $\ad(X)^3=-\ad(X)$ and an
% isometry $f:\R^n\to\R^d\times\frakp$ with $f(M)=\R^d\times\Ad(\rmK)\,X$\;.
% \end{theorem}

\subsection{Symmetric submanifolds associated with irreducible symmetric
  R-spaces}
\label{s:Naitoh} 
Continuing with the notation from the last section, we now
introduce certain symmetric submanifolds of $N$ which were
already mentioned in Theorem~\ref{th:Naitoh}. 
% Following~\cite{BENT}, if $N$ admits a symmetric R-space as described in
% Definition~\ref{de:symmetric_R-space}, then we associate a family of symmetric
% submanifolds of $N$ therewith; according to Theorem~\ref{th:Naitoh}, these families play an important role for the classification of
% symmetric submanifolds in ambient symmetric spaces of higher rank.
% let $\Gamma:T_oN\to\frakp$ denote its inverse (sometimes called the ``transvection
% map'' of $N$ at $o$).
In the following, note that the linear map $\pi_1:\fraki(N)\to
T_oN,X\mapsto X^*(o)$ is surjective and that we have
\begin{equation}\label{eq:frakk}
\frakk=\Menge{Y\in\fraki(N)}{\pi_1(Y)=0}; \end{equation}
hence $\pi_1|\frakp$ induces a linear isomorphism $\frakp\cong T_oN$\,.

\bigskip
\begin{proposition}\label{p:BENT}
Let $N$ be a symmetric space of compact type or of non-compact type and suppose
that there exists $X\in\frakp$ with $X\neq 0$ and
$\ad(X)^3=\epsilon\,\ad(X)$\,. 
% For each $c\in\R$ let $t\in]0,\pi[$ be uniquely determined by
%$-\cot(t)=c$ and consider the inner automorphism $\rho_t:\Iso(N)^0\to
%\Iso(N)^0$ given by $\rho_t(g):=\exp(t\,X)g\exp(-t\,X)$ for each $g\in\Iso(N)^0$\,. 
%Then the orbit $M:=M_c:=\rho_t(K)\,o$ 
Then there exists a family of symmetric submanifolds $M_c\subset N$\,, where $c$
ranges over $\R$\,, uniquely determined by the following properties:

For each $c\in\R$ we have $o\in M_c$\,, $T_oM_c=\pi_1(\frakp_\epsilon)$\,, $\bot_oM_c=\pi_1(\frakp_0)$\;,
and the second fundamental form $h$ of $M_c$ is characterized by
\begin{equation}\label{eq:fetth_2}
\forall Y\in\frakp_\epsilon:\fetth(\pi_1(Y))=c\,\pi_2(J\,Y)\;,
\end{equation}
where $J$ is the linear map from Lemma~\ref{le:ad(X)}\,.
Therefore, $M_0$ is totally geodesic in $N$\,, whereas for $c\neq 0$ the
submanifold $M_c$ is 1-full in $N$ with non-degenerate second fundamental form at $o$\,.
\end{proposition}

\begin{proof}
In case $N$ is of non-compact type, the existence of the family $\{M_c\}$ is established
in Theorem~2.3 of~\cite{BENT} (likewise, cf.~\cite{BCO},~Prop..\,9.3.8). The compact case can be handled by
similar methods, cf.~\cite{BCO}, Ch.~9.3. Furthermore, it is well known that a symmetric submanifold is
uniquely determined by its ``2-jet'' $(T_pM,h_p)$ at one point $p\in M$\,,
cf.~\cite{St}; thus $M_c$ is uniquely determined. For the last conclusion of the above theorem, notice that, as a result of
Proposition~\ref{p:symmetric_R-space}, both $M_c$ and the symmetric R-space
$\Ad(\rmK)\,X$ have the same tangent space at $o$ and at $X$\,, respectively (by means of identification $\frakp\cong T_oN$ via $\pi_1$).
Using also that $\pi_1$ is an equivariant map
of $\frakk$-modules, i.e.
\begin{equation}\label{eq:correspondence} 
\forall X\in\frakk,Y\in\frakp:\pi_1(\ad(X)\,Y)=\pi_2(X)\,\pi_1(Y)\;,
\end{equation} 
we now see from comparing~\eqref{eq:fetth_1} with~\eqref{eq:fetth_2} that
(again by means of the identification $T_oN\cong\frakp$) the
second fundamental form of $M_c$ at $o$ equals $c$ times the second
fundamental form of $\Ad(K)\,X$ at $X$\,; hence the result follows, again as a consequence of Proposition~\ref{p:symmetric_R-space}.
\end{proof}

\bigskip
\begin{definition}\label{de:BENT}
In the situation of Definition~\ref{de:symmetric_R-space}, 
suppose that there exists $X\in\frakp$ with $X\neq 0$ and $\ad(X)^3=\epsilon\,\ad(X)$\,. 
Then $N$ admits a symmetric R-space, 
and the family of submanifolds $M_c\subset N$ from Proposition~\ref{p:BENT} will be called  ``the family of
symmetric submanifolds associated therewith''.
\end{definition}

\section{Some intrinsic properties of $\osc f$}
\label{s:osc}
{\bf Throughout this section, $\mathbf{f:M\to N}$ is a parallel isometric immersion.}  
Then for each $p\in M$ the linear map $T_pf:T_pM\to
\osc_pf, x\mapsto T_pf\, x$ 
induces an injective vector bundle homomorphism
$TM\hookrightarrow f^*TN$ (whose image $Tf(TM)$ is usually called the ``first
osculating bundle of $f$''); hence we have the corresponding orthogonal splitting
\begin{equation}\label{eq:results-1}
f^*TN=Tf(TM)\oplus\bot f\cong TM\oplus\bot f\,. 
\end{equation}
In the following, in order to keep to our convention that ``$T_pM$ is
seen as a linear subspace of $T_{f(p)}N$'' for each $p\in M$\,, we
suppress the vector bundle isomorphism $Tf:TM\to Tf(TM)$\,;
for convenience, the reader may assume that $M\subset N$ is a submanifold and $f=\iota^M$\,.

\bigskip
\begin{definition}[Split-parallelity]
\label{de:split-parallel}
\begin{enumerate}
\item
The \emph{split connection} is by definition 
the linear connection  $\nabla^{sp}:=\nabla^M\oplus\nabla^\bot$\ on $f^*TN=TM\oplus\bot^1f$\,.
A section of $f^*TN$ will be called \emph{split-parallel} if it is parallel with respect to $\nabla^{sp}$\;.
\item
For a curve $c:J\to M$ 
let $\ghdisp{t_1}{t_2}{c}{sp}:T_{c(t_1)}N\to T_{c(t_2)}N$ denote the corresponding {\em split parallel displacement}\/ along $c$ (where $(t_1,t_2)$ varies
over $J\times J$), which 
is the family of linear isometries characterized by the following properties:
\begin{itemize}
\item $\ghdisp{t_1}{t_2}{c}{sp}\,X(t_1)=X(t_2)$ for any $\nabla^M$-parallel
  section $X:J\to TM$ along $c$\,,
\item $\ghdisp{t_1}{t_2}{c}{sp}\,\xi(t_1)=\xi(t_2)$ for any \,$\nabla^\bot$\,-parallel section $\xi:J\to \bot M$ along $c$\,.
\end{itemize}
\end{enumerate}
\end{definition}

Now the equations of Gau{\ss} and Weingarten~\eqref{eq:GW}
can formally be combined to
\begin{align}
&\forall\,X\in\Gamma(TM),V\in\Gamma(f^*(TN)):
\label{eq:first_Gauss}
\Kovv{N}{X}{(V)}=\Kovv{sp}{X}{V}+\fetth(X)\,V\;.
\end{align}
$R^{sp}$ and $R^\bot$ will denote the curvature tensors of $f^*TN$ and
$\bot f$ with respect to $\nabla^{sp}$ and $\nabla^\bot$\,, respectively. 
Since $f$ is parallel, also the curvature equations of Gau{\ss}, Codazzi and Ricci
can formally be combined to
\begin{equation}\label{eq:Gauss_Ricci}
\forall x,y\in T_pM:R^N(x,y)=R^{sp}(x,y) +[\fetth(x),\fetth(y)]\;.
\end{equation}
The following lemma is proved in a straightforward manner:

\bigskip
\begin{lemma}\label{le:splitting}
Let $V$ be a Euclidian vector space and $W\subset V$ a linear subspace.
Recall the splitting $\so(V)=\so(V)_+\oplus\so(V)_-$ defined by~\eqref{eq:even} and~\eqref{eq:odd}.
\begin{enumerate}
\item
We have $A\in\so(V)_+$ if and only if $A(W)\subset W$\,. 
\item
The map $\so(V)_-\to\rmL(W,W^\bot),A\mapsto A|W$ is a linear isomorphism.
\end{enumerate}
\end{lemma}

\bigskip
\begin{proposition}\label{p:locally_symmetric}
%\Abstand{2}
Let a parallel isometric immersion $f:M\to N$ be given.
\begin{enumerate}
\item $T_pM$  is a curvature invariant subspace of $T_{f(p)}N$\,, and we have 
\begin{equation}\label{eq:ci}
\forall x,y\in T_pM:\,R^N(x,y)\in\so(T_oN)_+\;.
\end{equation}
\item $M$ is locally symmetric, i.e. $R^M$ is parallel.
\item If $M$ is a complete, simply
connected, parallel submanifold, then $M$ is a symmetric space.
\item $h$ satisfies a second order tensorial property known as
  ``semiparallelity'':
\begin{align}\label{eq:semiparallel_1}
\forall x_1,x_2,y_1,y_2\in T_pM:\,R^\bot(x_1,x_2)\,h(y_1,y_2)=h(R^M(x_1,x_2)\,y_1,y_2)+h(y_1,R^M(x_1,x_2)\,y_2)\;.
\end{align}
\item Equation~\eqref{eq:semiparallel_1} is equivalent to
\begin{equation}\label{eq:semiparallel_2}
\forall x,y,z\in T_pM:\,\fetth(R^M(x,y)\,z)=[R^N(x,y)-[\fetth(x),\fetth(y)],\fetth(z)]\;.
\end{equation}
\end{enumerate}
\end{proposition}
\begin{proof}
(a) follows from the Codazzi Equation and Lemma~\ref{le:splitting}.
For the proof of (b) one needs assertion (a) and the curvature equation of
Gau{\ss}.
If $M$ is simply connected and complete, then it is even globally symmetric (cf.~\cite{He},~Ch.\,IV,~\S~6,~Theorem~5.6). 
The proof of Equation~\eqref{eq:semiparallel_1} is straightforward, see for example~\cite{Ts} or~\cite{Fe}.
For~\eqref{eq:semiparallel_2}, note that both sides of this equation are
elements of $\so(T_pN)_-$\,.
Thus by virtue of Lemma~\ref{le:splitting} it is enough to verify that
\eqref{eq:semiparallel_2} holds on $T_pM$\,. For this let $\tilde y\in T_pM$ be
given; then~\eqref{eq:Gauss_Ricci} implies:
$$
[R^N(x,y)-[\fetth(x),\fetth(y)],\fetth(z)]\,\tilde
y=R^\bot(x,y)\,\fetth(z)\,\tilde y-\fetth(z)\,R^M(x,y)\,\tilde y\;;
$$
now use~\eqref{eq:semiparallel_1}\;.
\end{proof}

\bigskip
\begin{proposition}\label{p:properties} 
Let a parallel isometric immersion $f:M\to N$ be given.
\begin{enumerate}
\item $\bot^1f\subset \bot f$ is a $\nabla^\bot$-parallel vector
subbundle. In particular,
\begin{equation}
\label{eq:R_bot}
\forall\;x,y\in T_pM:\;R^\bot(x,y)(\bot^1_pf)\subset\bot^1_pf\;.
\end{equation}
\item $\osc f$ is a
split-parallel vector subbundle of $f^*TN$\,.
\item If $c:J\to M$ is a curve, $X(t),Y(t)$ are
parallel sections of $TM$ and $\xi(t)$ is a parallel section of $\bot^1f$ along $c$\,,
then $R^\bot(X(t),Y(t))\,\xi(t)$ is a parallel section of $\bot^1f$ along $c$\,.
\item
$\osc f$ is a $\nabla^N$-parallel vector subbundle of $f^*TN$\,. Hence $\nabla^N$ induces a
connection on $\osc f$\,, as already described in Section~\ref{s:overview}.
%and a split-parallel
Therefore
\begin{align}\label{eq:osc_is_parallel_1}
&\forall t_1,t_2\in
J:\ghdisp{t_1}{t_2}{c}{N}(\osc_{c(t_1)}f)=\osc_{c(t_2)}f\;,\\
\label{eq:osc_is_parallel}
&\forall\;x,y\in T_pM: R^N(x,y)(\osc_pf)\subset\osc_pf\;;
\intertext{and the corresponding parallel displacement resp.\ curvature
  tensor are given by}
\label{eq:parisp_osc}
&\ghdisp{t_1}{t_2}{c}{\osc
f}=\big(\ghdisp{t_1}{t_2}{c}{N}\big)^{\osc}\qmq{(see~\eqref{eq:A_osc}),}\\
\label{eq:R_osc_1}
&R^{\osc f}(x,y)=\big(R^N(x,y)\big)^{\osc}\;.\\
\intertext{Moreover, we have}
\label{eq:R_osc_2}
& R^{\osc f}(x,y)\in\so(\osc_of)_+\;.
\end{align}
\item 
If $c:J\to M$ is a curve, $X(t),Y(t)$ are split-parallel sections of $TM$ and
$V(t)$ is a split-parallel section of $\osc f$ along $c$\,,
then $R^N(X(t),Y(t))V(t)$ is a split-parallel section of $\osc f$ along $c$\,, too\;.
\end{enumerate}
\end{proposition}
\begin{proof}
(a) and~(b) are straightforward,
(c) is proved using Part~(b) of
Proposition~\ref{p:locally_symmetric},~\eqref{eq:semiparallel_1} and~\eqref{eq:R_bot}.
(d) is again straightforward (for~\eqref{eq:R_osc_2} use the curvature
invariance of $T_pM$\,, Lemma~\ref{le:splitting} and~\eqref{eq:split_osc}). (e) follows from Part~(b) of
Proposition~\ref{p:locally_symmetric} and (c),
applying two times the curvature equations of Gau{\ss} and Ricci,
cf. Lemma~2.3 in~\cite{Ts}.
\end{proof}
% In order to avoid any cumbersome notation in the following, from now on we assume that $M\subset N$ is a parallel
% submanifold and $f=\iota^M$\,.

\bigskip
\begin{lemma}
In the situation of Definition~\ref{de:split-parallel}, we have for all $x,y\in T_{c(t_1)}M$\,:
\begin{align}
%\label{eq:osc}\\
\label{eq:h_is_parallel}
&h(\ghdisp{t_1}{t_2}{c}{M}\,x,\ghdisp{t_1}{t_2}{c}{M}\,y)=\ghdisp{t_1}{t_2}{c}{sp}\,h(x,y)\;,\\
\label{eq:fetth_is_parallel}
&\fetth(\ghdisp{t_1}{t_2}{c}{M}x)=\ghdisp{t_1}{t_2}{c}{sp}\circ\fetth(x)\circ{\ghdisp{t_2}{t_1}{c}{sp}}\qmq{and}\\
& \label{eq:oscM_is_split-parallel}
\forall t\in J:\ghdisp{t_1}{t_2}{c}{sp}(\osc_{c(t_1)}f)=\osc_{c(t_2)}f\;.
\end{align}
\end{lemma}
\begin{proof}
\eqref{eq:h_is_parallel} follows by definition of the parallelity of the second fundamental form.
In particular, \eqref{eq:fetth_is_parallel} holds on $T_{c(t_1)}M$\,; moreover both sides of Equation~\eqref{eq:fetth_is_parallel} are
elements of $\so(T_{c(t)}N)_-$ (because of Equation~\eqref{eq:fetth_in_osc}
and since $\ghdisp{t_1}{t_2}{c}{sp}$ respects the two splittings
$T_{f(c(t_1))}N=T_{c(t_1)}M\oplus\bot_{c(t_1)}M$ and $T_{f(c(t_2))}N=T_{c(t_2)}M\oplus\bot_{c(t_2)}f$), and thus~\eqref{eq:fetth_is_parallel} follows from Lemma~\ref{le:splitting}.
%As $\fetth$ and $h$ can be canonically identified using only the metric and the splitting $TN|M=TM\oplus \bot M$\,, (b) also follows from (a).
Equation~\eqref{eq:oscM_is_split-parallel} follows immediately from Part~(b) of Proposition~\ref{p:properties}.
\end{proof}

Because of~\eqref{eq:fetth_is_parallel}, the following proposition is a
consequence of Lemma~3 in~\cite{JR}:

\bigskip
\begin{proposition}\label{P:Move}
For each curve $c:J\to M$ with $0\in J$ and $c(0)=p$ let $X$ denote the backward parallel transport of the velocity vector field
$\dot c$\,, i.e.
\begin{equation}\label{eq:X(t)}
X(t):=\ghdisp{t}{0}{c}{M}(\dot c(t))\in T_p M\;.
\end{equation}

Then the function
\begin{equation}\label{Eq:Move2}
\mu_c:J\to \SO(T_{c(0)}N) \ ,\ t\mapsto\ghdisp{t}{0}{c}{N}\circ\ghdisp{0}{t}{c}{sp}
\end{equation}

solves the linear differential equation
\begin{equation}\label{Eq:Move3}
\mu_c'(t)\:=\mu_c(t)\circ \fetth(X(t)) \qmq{with}\mu_c(0)=\id\;.
\end{equation}
\end{proposition}

Equations~\eqref{eq:oscM_is_split-parallel} and~\eqref{eq:osc_is_parallel_1} imply
for a curve $c:[0,1]\to M$ as above:
\begin{equation}
\label{eq:mu_c_1}
\forall t\in J:\mu_c(t)(\osc_pf)=\osc_pf\;.
\end{equation}
Moreover, using the canonical identification $$\SO(\osc_pf)\cong\Menge{g\in
  \SO(T_pM)}{g(\osc_pf)=\osc_pf\qmq{and}g|(\osc_pf)^\bot=\Id}\;,$$ 
by means of~\eqref{eq:fetth_in_osc} and~\eqref{Eq:Move3} we have
\begin{equation}
\label{eq:mu_c_2}
\forall t\in J:\mu_c(t)\in\SO(\osc_pf)
\end{equation}

\bigskip
\begin{example}\label{ex:mu}
\begin{enumerate}
\item
If $c$ denotes the geodesic $\gamma_x:J\to M$ with $\dot \gamma_x(0)=x$\,,
then
\begin{equation}\label{eq:mu_1}
\mu_{c}(t)=\exp(t\,\fetth(x))\;.
\end{equation}
\item
We can construct $\mu_c$ as in Proposition~\ref{P:Move} also if we only assume
that $c$ is continuous and piecewise differentiable; then $\mu_c$ will
be continuous and piecewise differentiable, too. Now let $x,y\in T_pM$ be
given, consider the
corresponding smooth geodesic line $\gamma_x:\R\to M$ and put $$\tilde
y:=\ghdisp{0}{1}{\gamma_x}{M}\,y\in T_{\gamma_x(1)}M\;;$$
 hence we also have the
corresponding smooth geodesic line $\gamma_{\tilde
  y}:\R\to M$\,.
Let $\gamma_{(x,y)}:[0,2]\to M$  be the broken geodesic line characterized by
\begin{equation}\label{eq:gamma_xy}
\forall t\in[0,2]:\gamma_{(x,y)}(t)=
\begin{cases}
\gamma_x(t)&\text{\ for }t\leq 1\\
\gamma_{\tilde y}(t-1)&\text{\ for }t\geq 1
\end{cases}
\end{equation}
Then for the curve given by $c(t):=\gamma_{(x,y)}(t)$ we have
 \begin{equation}\label{eq:mu_2}
\mu_{c}(2)=\exp(\fetth(x))\circ\exp\big(\fetth(y)\big)\;.
\end{equation}
\end{enumerate}
\end{example}
\begin{proof}
For (a): Here the function $X(t)$ of Equation~\eqref{eq:X(t)} is constant
equal to $x$\,;
therefore the solution of~\eqref{Eq:Move3} is the one-parameter subgroup given by~\eqref{eq:mu_1}.

For (b): 
By virtue of Equation~\eqref{Eq:Move2} and
using Part~(a) several times, we obtain 
\begin{align*}
\mu_c(2)&=\ghdisp{2}{0}{c}{N}\circ\ghdisp{0}{2}{c}{sp}
=\ghdisp{1}{0}{\gamma_x}{N}\circ\ghdisp{1}{0}{\gamma_{\tilde y}}{N}\circ\ghdisp{0}{1}{\gamma_{\tilde
  y}}{sp}\circ\ghdisp{0}{1}{\gamma_x}{sp}\\
&=\ghdisp{1}{0}{\gamma_x}{N}\circ\ghdisp{0}{1}{\gamma_x}{sp}\circ\ghdisp{1}{0}{\gamma_x}{sp}\circ
\ghdisp{1}{0}{\gamma_{\tilde y}}{N}\circ\ghdisp{0}{1}{\gamma_{\tilde y}}{sp}\circ\ghdisp{0}{1}{\gamma_{x}}{sp}
\\
&=\exp(\fetth(x))\circ\exp\big(\ghdisp{1}{0}{\gamma_x}{sp}\circ \fetth(\tilde y)\circ\ghdisp{0}{1}{\gamma_x}{sp}\big)\;.
\end{align*}
The result follows from the previous together with~\eqref{eq:fetth_is_parallel}\,.
\end{proof}

\bigskip
\begin{lemma}\label{le:R(t)}
In the situation of Proposition~\ref{P:Move}, for any choice of
vectors $y_1,y_2\in T_pM$ and $v\in\osc_pf$ we have
\begin{align}\notag
&R^N(\mu_c(t)\,y_1,\mu_c(t)\,y_2)\,v\in\osc_pf\,,\\
\intertext{and the following two equalities hold:}
\label{eq:R(t)_3}
&R^N(\mu_c(t)\,y_1,\mu_c(t)\,y_2)\,v=\ghdisp{t}{0}{c}{N}\circ
R^N(\ghdisp{0}{t}{c}{sp}\,y_1,\ghdisp{0}{t}{c}{sp}\,y_2)\circ
\ghdisp{0}{t}{c}{N}\,v\;,\\
&R^N(\mu_c(t)\,y_1,\mu_c(t)\,y_2)\,v=\mu_c(t)\circ R^{N}(\,y_1,\,y_2)\circ\mu_c(t)^{-1}\,v\;.
    \label{eq:R(t)_2}
\end{align}
\end{lemma}
\begin{proof}
Using the $\nabla^N$-parallelity of $R^N$ and Part~(e) of
Proposition~\ref{p:properties}, we have
\begin{align*}\notag
&R^N(\mu_c(t)\,y_1,\mu_c(t)\,y_2)\,\mu_c(t)\,v
=\ghdisp{t}{0}{c}{N}\,R^N(\ghdisp{0}{t}{c}{sp}\,y_1,\ghdisp{0}{t}{c}{sp}\,y_2)\,\ghdisp{0}{t}{c}{sp}\,v\\
&=\ghdisp{t}{0}{c}{N}\ghdisp{0}{t}{c}{sp}\,R^N(y_1,y_2)\,v=\mu_c(t)(R^N(y_1,y_2)\,v)\;.
\end{align*}
The result follows with~\eqref{eq:mu_c_1}\;.
\end{proof}

\bigskip
For every Euclidian vector space $V$ and every $A\in \so(V)$ let
$A^{(2)}:\wedge^2(V)\to\wedge^2(V)$ 
denote the induced Endomorphism, i.e.
\begin{equation}\label{eq:derive_1}
%A\cdot R:=-\frac{\diff}{\diff t}|_{t=0}\exp(t\,A)^*R\;.
A^{(2)}\,u\wedge v:=A\,u\wedge v+u\wedge A\,v
=\frac{\diff}{\diff t}\Big|_{t=0}\exp(t\,A)\,u\wedge \exp(t\,A)\,v\;.
\end{equation}
Furthermore, one knows that for $A,B\in \so(V)$
\begin{equation}\label{eq:derive_3}
\frac{\diff}{\diff t}\Big|_{t=0}\exp(t\,A)\circ B\circ\exp(t\,A)^{-1}=[A,B]\;.
\end{equation}

\bigskip
\begin{lemma}\label{le:Lemma}
For arbitrary $p\in M$\,, $x_1,x_2\in T_pM$\,, $y_1,y_2\in T_pM$ and $v\in\osc_pf$
we have:
\begin{align}\label{eq:parallel_1}
&[\fetth(x_1),R^{N}(y_1,y_2)]\,v=R^N(\fetth(x_1)^{(2)}\,y_1\wedge y_2)\,v\;,\\
\label{eq:parallel_2}
&[\fetth(x_1),[\fetth(x_2),R^{N}(y_1,y_2)]]\,v
=R^N(\fetth(x_1)^{(2)}\fetth(x_2)^{(2)}\,y_1\wedge y_2)\,v\;.
\end{align}
\end{lemma}
\begin{remark}
Suppose that $\fraki(N)$ is semisimple and let $\frakk$ denote the Lie algebra
of the isotropy group at $p$\,. In accordance with~\cite{He},~Ch.\,V,~Theorem~4.1., we have
\begin{equation}\label{eq:derive}
\pi_2(\frakk)=\Menge{A\in\so(T_pN)}{\forall u,v\in T_pN:[A,R^N(u,v)]=R^N(Au,v)+R^N(u,Av)}\;.
\end{equation}
Hence, if $M$ is a symmetric submanifold of $N$\,, then, according to~\eqref{eq:Eschenburg}, Equation~\eqref{eq:parallel_1} (and
therefore also~\eqref{eq:parallel_2}) holds for \emph{all} $y_1,y_2\in T_pN$ \,.
\end{remark}
\begin{proof}[Proof of Lemma~\ref{le:Lemma}]
To derive Equation~\eqref{eq:parallel_1},
let $\gamma:J\to M$ denote the geodesic with $\dot\gamma(0)=x_1$\,.
Then $\mu_\gamma(t)=\exp(t\,\fetth(x_1))$ holds by means of Equation~\eqref{eq:mu_1};
considering also Equations~\eqref{eq:derive_1} and \eqref{eq:derive_3} (with
$A=\fetth(x_1)$)~\eqref{eq:parallel_1} therefore follows
by taking the derivative $\frac{\diff}{\diff t}\big|_{t=0}$ on both sides of
Equation~\eqref{eq:R(t)_2}, 

To derive Equation~\eqref{eq:parallel_2}\footnote{It should be mentioned that this equation can not be obtained from
  the previous one by iteration,
because $\fetth(x_2)^{(2)}\,y_1\wedge y_2$ is not an element of $\Lambda^2T_pM$\,.},
let $(s,r)\in \R\times \R$ be fixed elements and $c$ be the broken geodesic
$
\gamma_{(s\,x_1,r\,x_2)}:[0,2]\to M\,,
$
as described in Equation~\eqref{eq:gamma_xy}. Then $\mu_c$ satisfies
$\mu_c(2)=\exp(s\,\fetth(x_1))\circ \exp(r\,\fetth(x_2))=:f(s,r)$\,, according
to~\eqref{eq:mu_2}\,;
hence Equation~\eqref{eq:R(t)_2} gives
$$
\big(R^N(f(s,r)\,y_1,f(s,r)\,y_2)\big)\,v=\big(f(s,r)\circ R^N(y_1,y_2)\circ f(s,r)^{-1}\big)\,v\;.
$$
Now \eqref{eq:parallel_2} follows by taking the derivatives $\frac{\del}{\del
  r}\frac{\del}{\del s}\big|_{r=s=0}$
on both sides of this equation.
\end{proof}

\subsection{Curvature invariance of the first normal spaces}
\begin{proposition}\label{p:fundamental}
For arbitrary $p\in M$\,, $x,y\in T_pM$ we have
$\fetth(x)\,\fetth(y)(T_pM)\subset T_pM$ and the following equation holds on $\osc_pf$\,:
\begin{align}
%\notag
&R^N(h(x,x),h(y,y))
\label{eq:fundamental}
=[\fetth(x),[\fetth(y),R^N(x,y)]]
-R^N(\fetth(x)\,\fetth(y)\,x,y)-R^N(x,\fetth(x)\,\fetth(y)\,y)\;.
\end{align}
Moreover, for all $\xi,\eta\in\bot^1_pf$ the curvature endomorphism
$$
R^N(\xi,\eta):T_pN\to T_pN, v\mapsto R^N(\xi,\eta)\,v
$$
has the following property:
\begin{align}\label{eq:fundamental_2}
&R^N(\xi,\eta)(\osc_pf)\subset\osc_pf\qmq{and}\\
\label{eq:fundamental_1}
&\big(R^N(\xi,\eta)\big)^{\osc}\in \so(\osc_pf)_+\;.
\end{align}
\end{proposition}
\begin{proof}
Let us first verify Equation~\eqref{eq:fundamental} on $\osc_pf$\,:
According to~\eqref{eq:parallel_2}, we have
$$
\forall v\in\osc_pf:\;[\fetth(x),[\fetth(y),R^N(x,y)]]\,v=R^N(\fetth(x)^{(2)}\,\fetth(y)^{(2)}\,x\wedge y)\,v\;,
$$
and furthermore (using~\eqref{eq:derive_1} twice and the symmetry $h(y,x)=h(x,y)$):
\begin{align*}
&\fetth(x)^{(2)}\fetth(y)^{(2)}\,x\wedge y=\fetth(x)\,\fetth(y)\,x\wedge y
+\underbrace{\fetth(x)\,x\wedge\fetth(y)\,y}_{=h(x,x)\wedge h(y,y)}
+\underbrace{\fetth(y)\,x\wedge\fetth(x)\,y}_{=h(y,x)\wedge h(x,y)=0} +x\wedge\fetth(x)\,\fetth(y)\,y\;.
\end{align*} Therefore, Equation~\eqref{eq:fundamental} holds on $\osc_pf$\,. 
For the proof of \eqref{eq:fundamental_2}, it is enough to assume that there
exist $x,y\in T_pM$
with $\xi=h(x,x)$\,, $\eta=h(y,y)$\,, because $h$ is a symmetric bilinear
map. Furthermore, note that
$\fetth(x)\,\fetth(y)\,z=-S_{h(y,z)}x$ for all $z\in T_pM$\,, hence
$\fetth(x)\,\fetth(y)(T_pM)\subset T_pM$ and therefore the linear space $\osc_pf$ is
invariant under each of the three terms on the right hand side
of~\eqref{eq:fundamental}, in accordance with~\eqref{eq:so_osc},~\eqref{eq:fetth_in_osc},~\eqref{eq:osc_is_parallel}; which implies Equation~\eqref{eq:fundamental_2}\;.
To conclude also~\eqref{eq:fundamental_1}, just note that (after
projection to $\so(\osc_of)$) each of the
three terms on the right hand side
of~\eqref{eq:fundamental} is an element of $\so(\osc_pf)_+$\,, according to~\eqref{eq:fetth_in_osc},~\eqref{eq:R_osc_1} and~\eqref{eq:R_osc_2}
and the rules for $\Z_2$-graded Lie algebras.
\end{proof}
As a consequence of~\eqref{eq:even},~\eqref{eq:fundamental_2} and~\eqref{eq:fundamental_1} we have:

\bigskip
\begin{corollary}\label{co:ci}
$\bot^1_pf$ is a curvature invariant subspace of $T_{f(p)}N$\,.
\end{corollary}
Motivated by Lemma~2 in~\cite{E}, we are now able to generalize Part~(e) of Proposition~\ref{p:properties}:

\bigskip
\begin{proposition}\label{p:R_is_split-parallel}
If $f:M\to N$ is parallel, then the tensor of type $(0,4)$ on $\osc f$ defined by
\begin{equation}\label{eq:R|osc}
R^\flat(v_1,v_2,v_3,v_4):=\g{R^N(v_1,v_2)v_3}{v_4}\qmq{for} v_1,\ldots,v_4\in \osc_pf
\end{equation}
is  split-parallel; which means: For a curve $c:J\to M$ and split-parallel
sections of $\osc f$ along $c$\,, $V_i(t)$
($i=1,\ldots,4$), the function $f(t):=R^\flat\big(V_1(t),V_2(t),V_3(t),V_4(t)\big )$ is constant\;.
\end{proposition}
\begin{proof}
Note that $V_i(t)=X_i(t)+\xi_i(t)$\,, where $X_i$ resp.\ $\xi_i$ are
split-parallel sections of $TM$ resp.\ of $\bot^1f$\,.
%For the proof put $f(t):=R^\flat\big(V_1(t),V_2(t),V_3(t),V_4(t)\big )$\,.
Thus it is enough to consider the following two cases.

First case: Exactly one of the sections $V_i$
%, say $V_1$ (otherwise use the symmetries of $R^N$),
is a section of $\bot^1f$ (resp.\ of $TM$) and the other ones are sections of $TM$
(resp.\ of $\bot^1f$).
Then $f(t)=0$\,, since $T_{c(t)}M$ (resp.\ $\bot^1_{c(t)}f$) are
curvature invariant subspaces of $T_{c(t)}N$\;;
see Proposition~\ref{p:locally_symmetric} (a) and the previous Lemma.
Note that for this argument it was not used that the sections are split-parallel.

Second case: An even number of $V_1,V_2,V_3,V_4$ are sections of $TM$\,,
and the other ones are sections of $\bot^1f$\,. Then, by the equations of Gau{\ss}
and Weingarten, for each $i\in \{1,\ldots,4\}$ an odd number of the sections
$V_1,\ldots,\Kodel{N}{V_i},\ldots V_4$ are sections of $TM$\,, and the other ones
are sections of $\bot^1f$\,.
It follows from the parallelity of $R^N$ and the considerations made for the first case that $f'(t)=0$\,.
\end{proof}

\bigskip
\begin{lemma}\label{le:more_general}
Equations~\eqref{eq:R(t)_3},~\eqref{eq:R(t)_2},~\eqref{eq:parallel_1}
and~\eqref{eq:parallel_2} also hold if one replaces ``$y_1,y_2\in T_pM$'' with
``$y_1,y_2\in\osc_pf$'' in each of these Equations.
\end{lemma}
\begin{proof}
The result follows by repeating the proofs of~\eqref{eq:R(t)_3},~\eqref{eq:R(t)_2},~\eqref{eq:parallel_1}
and~\eqref{eq:parallel_2}, but now using Proposition~\ref{p:R_is_split-parallel} instead of
Proposition~\ref{p:properties}~(e).
\end{proof}
Lemma~\ref{le:more_general} (applied to~\eqref{eq:parallel_2}) will be
needed for the proof of Theorem~\ref{th:hol} in Section~\ref{s:hol}.

\bigskip
\begin{corollary}
%$\fetth(x)$ is a derivation of $R^\flat$\,.
For all $x\in T_pM$ and $v_1,\ldots,v_4\in \osc_pf$ we have
\begin{equation}\label{eq:fetth_3}
\sum_{i=1}^4R^{\flat}(v_1,\ldots,\fetth(x)\,v_i,\ldots,v_4)=0\;.
\end{equation}
\end{corollary}
\begin{proof}
From Lemma~\ref{le:more_general} (applied to Equation~\eqref{eq:R(t)_2}) we obtain
\begin{equation}\label{eq:R_b_is_constant}
R^{\flat}(\mu_c(t)\,v_1,\mu_c(t)\,v_2,\mu_c(t)\,v_3,\mu_c(t)\,v_4)=const
\end{equation}
for each curve $c:\R\to M$\,. If $c=\gamma_x$ is the geodesic considered in
Example~\ref{ex:mu}, then~\eqref{eq:fetth_3} follows with~\eqref{eq:mu_1}
by taking the derivative $\frac{\diff}{\diff t}\big|_{t=0}$ of~\eqref{eq:R_b_is_constant}.
\end{proof}

\subsection{Proof of  Theorem~\ref{th:result_1}}
\label{s:1-full_versus_symmetric}
At the end of this section, we will give the proof of Theorem~\ref{th:result_1}. 

\begin{lemma}\label{le:non-degenerate}
Let a full parallel submanifold $M\subset N$ with $o\in M$ be given and $\hol(M^\top)$ be the holonomy Lie algebra of the totally geodesic
submanifold $M^\top:=\exp^N(T_oM)$ with respect to the base point $o$\,. Then $\Kern(\fetth^M_o)$ is a $\hol(M^\top)$-invariant subspace; hence,  
if $M^\top$ is an {\em irreducible}\/
symmetric space, then $h^M_o$ is a non-degenerate symmetric bilinear form.
\end{lemma}
\begin{proof}
According to the Theorem of Ambrose/Singer, since $M^\top$ is a totally geodesic submanifold of $N$\,, $\hol(M^\top)$ is the linear subspace of
$\so(T_pM)$ given by
$$
\Menge{R^N(x,y)|T_oM:T_oM\to T_oM}{x,y\in T_oM}_\R\;.
$$
I claim that 
$
\Kern(\fetth_o)
$
is a subspace of $T_oM$ which is invariant under the action of $\hol(M^\top)$\,; hence,
if $\hol(M^\top)$ acts irreducible on $T_oM$\,, then $\fetth$ is injective,
because $\fetth=0$ is not possible for a full submanifold.

For this, let $z\in T_oM$ be given and assume that $\fetth(z)=0$\,. I claim
that then also $\fetth(A(z))=0$ for each $A\in\hol(M^\top)$\,.

For this: By the previous, we may assume that there exist $x,y\in T_oM$
with $A=R^N(x,y)|T_oM:T_oM\to T_oM$\,. Then r.h.s.\ of~\eqref{eq:semiparallel_2}
  vanishes, therefore $\fetth(R^M(x,y)\,z)=0$ and thus we have (by the
  Gau{\ss} equation for the curvature)
 $$
\fetth(R^N(x,y)z)\stackrel{\eqref{eq:Gauss_Ricci},\eqref{eq:semiparallel_2}}
=\fetth\big(\fetth(x)\!\!\underbrace{h(y,z)}_{=h(z,y)=0}-\fetth(y)\underbrace{h(x,z)}_{=0}\big)=0\;.
$$
The result follows.
\end{proof}

\bigskip
\begin{lemma}\label{le:ci}
% \Abstand{3}
For the following types
of parallel submanifolds their second osculating spaces are curvature invariant:
\begin{enumerate}
\item $M$ is a submanifold of a real space form.
\item $M^m$ is a totally real submanifold of the complex projective space
  $\C\rmP^n$ and $m>1$\,, see~\cite{N1}, Lemma 2.1.
\item $M$ is a complex submanifold of $\C\rmP^n$ 
(since here the second osculating spaces are complex subspaces). 
%\,, see~\cite{Ts}.
% \item $M\subset N$ is a totally complex submanifold of 
%  (;
%   
\item $M$ is a totally complex submanifold of the quaternionic
projective space $\bbH\rmP^n$ which is locally of ``K\"{a}hlerian type'' (in
the sense of~\cite{Ts},~Definition~2.12), see~\cite{ADM},~Prop. 5.6\,.
\end{enumerate}
Furthermore, according to Proposition~2.11
of~\cite{Ts3}, a totally complex submanifold $M^{2m}\subset\bbH\rmP^n$ with $m\geq 2$ is already
locally of K\"{a}hlerian type.
\end{lemma}

\bigskip
\begin{proof}[Proof of Theorem~\ref{th:result_1}]
Part~(a) of Theorem~\ref{th:result_1} stands in accordance with
Corollary~\ref{co:ci}.

For Part~(b) and~(c): A 1-full, complete parallel submanifold of $N$ is even a symmetric submanifold as a
consequence of Part~(a) combined with Theorem~\ref{th:struebing}.

Conversely, let $M$ be a full symmetric submanifold of $N$\,; hence the
subgroup of $\Iso(N)$ generated by the extrinsic symmetries of $M$
acts transitively $M$ and therefore $M$ is a complete Riemannian
manifold. Moreover, its second fundamental form is parallel according to Theorem~\ref{th:struebing}.
It remains to show that $M$ is 1-full and, in case $N$ has no Euclidian
factor, that $h$ is non-degenerate: 

For this, because of Theorem~\ref{th:Naitoh_2}, it is
enough to assume that $N$ is an Euclidian space or an irreducible symmetric
space. If $N$ is a Euclidian space or if $N$ is irreducible and the rank of
$N$ is 1, then by virtue of Theorem~\ref{th:Naitoh} combined with
Lemma~\ref{le:ci}, the second osculating spaces of $M$ are curvature invariant; thus $M$ is even 1-full by means of
``reduction of the codimension'' (Theorem~\ref{th:Dombi}). 
Furthermore, if the rank of $N$ is 1, then for each $p\in M$ the totally
geodesic submanifold $M^\top(p)$ defined in Lemma~\ref{le:non-degenerate} is
either a real space form (in Cases~(a) and~(b) of Lemma~\ref{le:ci}) or a
complex space form (in Cases~(c) and~(d) of Lemma~\ref{le:ci}), and therefore
the non-degeneracy of $h_o$ is given by Lemma~\ref{le:non-degenerate}.
If $N$ is irreducible and the rank of $N$ is larger than 1, 
then by the strength of Theorem~\ref{th:Naitoh}, $N$ admits a symmetric
R-space and $M$ belongs to the family of symmetric submanifolds 
associated therewith, i.e. $M=M_c$ for some $c\in \R$\,. Then $c\neq
0$ (since $M$ is not totally geodesic in $N$) and hence $M$ is a 1-full
submanifold of $N$ and $h_o$ is non-degenerate according to
Proposition~\ref{p:BENT}.
\end{proof}

\section{Symmetry of $\osc f$}
\label{sec:symmetric_strip}
{\bf Throughout this section, $\mathbf{f:M\to N}$ is
  a parallel isometric immersion defined on a simply connected Riemannian
  symmetric space $\mathbf{M}$} 
(cf. Part (c) of Proposition~\ref{p:locally_symmetric}).
The corresponding geodesic symmetries of $M$ will be denoted by
$\sigma^M_p\;(p\in M)$\,. Remember that $\osc f\subset f^*TN$ is a
$\nabla^N$-parallel vector subbundle, according to Part~(d) of
Proposition~\ref{p:properties}, and hence $\osc f$ is equipped with the
connection $\nabla^{\osc f}$ induced by restriction of $\nabla^N$\,. Furthermore, continuing
with the notation from Section~\ref{s:osc}, there is the splitting
$\osc f=TM\oplus\bot^1f$\,; but note that in general $\nabla^{\osc f}$ is
not the split-connection (introduced in Definition~\ref{de:split-parallel})
restricted to $\osc f$\,. In Section~\ref{s:th_wes} we will prove the following result:

\bigskip
\begin{theorem}\label{th:wes}
For each $p\in M$ there exists 
a unique involutive map $\Sigma_p:\osc f\to \osc f$ characterized by the following properties:
% \Abstand{3}
\begin{enumerate} 
\item
  $\Sigma_p$ is a fibrewise isometric
vector bundle homomorphism along
$\sigma^M_p$\,, i.e.\ the following diagram is commutative
$$\begin{CD}
\osc f @>\Sigma_p>> \osc f\\
@VVV @VVV \\
M @>\sigma^M_p>> M\;
\end{CD}
$$
and for each $q\in M$ the map $\Sigma_p|\osc_qf:\osc_qf\to\osc_{\sigma^M_p(q)}f$ is a linear isometry.
\item $\Sigma_p$ is a $\nabla^{\osc f}$-parallel vector bundle isomorphism of $\osc f$\,.
\item
$\Sigma_p|\osc_pf$ is the linear reflection in $\bot^1_pf$\,.
\end{enumerate}
Moreover:
\begin{enumerate}\addtocounter{enumi}{3}
\item
\begin{equation}\label{eq:Sigma_restricted_to_TM}
\forall q\in M, x\in T_qM:\Sigma_p\,x=T\sigma_p^M\,x\;.
\end{equation}
\begin{equation}\label{eq:Sigma_is_split-parallel}
\forall q\in M,\forall x,y\in
T_qM:\Sigma_p\,h(x,y)=h(T\sigma_p^M\,x,T\sigma_p^M\,y)\;.
\end{equation}
\item
For every smooth geodesic line $\gamma$ of $M$ with $\gamma(0)=p$ we have $\sigma_p^M(\gamma(-1))=\gamma(1)$\,,
\begin{align}\label{eq:Sigma_2}
&\Sigma_p|T_{\gamma(-1)}M=-\ghdisp{-1}{1}{\gamma}{M}\;,\\
\label{eq:Sigma_3}
&\Sigma_p|\bot^1_{\gamma(-1)}f={\ \ \ }\ghdisp{-1}{1}{\gamma}{\bot}|\bot^1_{\gamma(-1)}f\;.
\end{align}
\item
If $M$ is a symmetric submanifold of $N$ (with extrinsic symmetries $\sigma_p^\bot$ ($p\in
M$)), then we have $T\sigma_p^\bot(\osc M)\subset\osc M$ and
$$
\Sigma_p|\osc M=T\sigma_p^\bot|\osc M:\osc M\to\osc M\;.
$$
\end{enumerate}
\end{theorem}
Comparing the last theorem with Definition~\ref{de:extrinsically_symmetric}, we
hence see that the family $\Sigma_p$ ($p\in M$) can be seen as sort of ``weak extrinsic symmetries'' of
$M$\,, and hence $M$ is (at least) ``extrinsically symmetric in $\osc f$''.

\bigskip
\begin{definition}\label{de:Phi}
For each geodesic $\gamma$ of $M$ with $\gamma(0)=p$ and each $t\in \R$ we define
\begin{equation}\label{eq:theta_and_Phi}
\theta_\gamma(t):=\sigma_{c(t/2)}^M\circ
\sigma_p^M\qmq{and}\Theta_\gamma(t):=\Sigma_{\gamma(t/2)}\circ \Sigma_p\;\;.
\end{equation}
\end{definition}
Then $\theta_\gamma(t)$ ($t\in\R$) is a family of isometries defined on $M$ and 
$\Theta_\gamma(t)$ ($t\in\R$) is a family of isometric, parallel vector bundle
isomorphism of $\osc f$ (by virtue of Theorem~\ref{th:wes}), and the following diagram is commutative:
\begin{equation}\label{eq:Phi}
\begin{CD}
\osc f @>\Theta_\gamma(t)>>\osc f \\
@VVV @VVV \\
M @>\theta_\gamma(t)>> M
\end{CD}
\end{equation}

\bigskip
\begin{corollary}\label{co:Phi}
In the situation of Definition~\ref{de:Phi}, 
let $\ghdisp{0}{t}{\gamma}{sp}$ denote the corresponding split-parallel
displacement 
along $\gamma$ as introduced in Definition~\ref{de:split-parallel}.
For each $t\in\R$ we have $\theta_\gamma(t)(p)=\gamma(t)$ and
\begin{align}
\label{eq:Phi_gleich_phi}
&\Theta_\gamma(t)|\osc_pf=\ghdisp{0}{t}{\gamma}{sp}|\osc_pf:\osc_pf\to \osc_{\gamma(t)}f\;.
\end{align}
\end{corollary}
\begin{proof}
Using Definition~\ref{de:Phi} and Part~(c) of Theorem~\ref{th:wes}, we have
$\theta_\gamma(t)(p)=\sigma_{\gamma(t/2)}^M(p)=\gamma(t)$ and for all
$x+\xi\in T_pM\oplus\bot^1_pf$
$$
\Theta_\gamma(t)(x+\xi)
=\Sigma_{\gamma(t/2)}(-x+\xi)
\stackrel{\eqref{eq:Sigma_2},\eqref{eq:Sigma_3}}=\ghdisp{0}{t}{\gamma}{sp}(x+\xi)\;,
$$
which yields the stated result.
\end{proof}

\subsection{Certain involutions on  the first normal bundle}
\label{s:th_wes}
At the end of this section we will give the proof of Theorem~\ref{th:wes}.
But first we have to verify the existence of certain involutions on
$\bot^1f$\,, for which purpose
we will now state some general facts about the existence of parallel sections of
some vector bundle $\bbE$ over a simply connected Riemannian manifold $M$  equipped with a connection.
Let $o\in M$ be a fixed ``origin'' and
$s_0\in\bbE_o$ considered as ``initial condition''.

\bigskip
\begin{lemma}\label{le:R_E_is_parallel}
Suppose that the curvature tensor $R^\bbE$ is parallel (considered as a section of $\rmL(\Lambda^2(TM),\End(\bbE))$\,, where the
latter space is equipped with the induced connection).
%$s_0$ is a fixed point with respect to the action of $\Hol(\bbE)$ if and only if
Then there exists a parallel section $s$ of $\bbE$ with $s(o)=s_0$ if and only if
\begin{equation}
    \forall x,y\in T_oM:\;R^\bbE(x,y)\,s_0=0\;.
    \label{eq:existence_2}
\end{equation}
\end{lemma}
\begin{proof}
Let $\Hol(\bbE)$ denote the holonomy group of $\bbE$ with respect to the base point $o$\,, defined by
$$
\Hol(\bbE):=\Menge{\ghdisp{0}{1}{c}{\bbE}}{c:[0,1]\to M\qmq{is a curve with}c(0)=c(1)=o}\;,$$
where $\ghdisp{0}{1}{c}{\bbE}$ means the parallel displacement along $c$ in $\bbE$\,.
It is known that $\Hol(\bbE)$ is a Lie subgroup of $\GL(\bbE_o)$\,,
its Lie algebra, $\hol(\bbE)\subset \End(\bbE_o)$\,, is called the holonomy
Lie algebra of $\bbE$\,. The ``Theorem of Ambrose/Singer'' shortly states that $\hol(\bbE)$
is generated by the curvature of $\bbE$\,; more exactly it is generated
(as a vector space over $\R$) by the elements $$\ghdisp{1}{0}{c}{\bbE}\circ
R^\bbE(x,y)\circ \ghdisp{0}{1}{c}{\bbE}\;,$$
where $c$ runs over all curves $[0,1]\to M$ with $c(0)=o$ and $x,y\in T_{c(1)}M$\,.
If $R^\bbE$ is
parallel, then we have
$$\ghdisp{1}{0}{c}{\bbE}\circ R^\bbE(x,y)\circ \ghdisp{0}{1}{c}{\bbE}=R^\bbE(\ghdisp{1}{0}{c}{M}\,x,\ghdisp{1}{0}{c}{M}\,y)\;;
$$
and therefore\begin{equation}\label{eq:hol_of_bbE}
\hol(\bbE)=\Spann{R^\bbE(x,y)}{x,y\in T_oM}\;.
\end{equation}
Let $s\in\Gamma(\bbE)$ be a section with $s(o)=s_0$\,. Then $s$
is parallel if and only if for every curve $c:[0,1]\to M$ with $c(0)=o$ we have
\begin{equation}\label{eq:s}
s(c(1))=\ghdisp{0}{1}{c}{\bbE}\,s_0\;.
\end{equation}
Thus, if there exists a parallel section with $s(o)=s_0$\,, then in particular $s_0$ is a fix point of
$\Hol(\bbE)$\,. And if $s_0$ is a fix point of
$\Hol(\bbE)$\,, then one defines a section $s$ via~\eqref{eq:s}, which is then
parallel with $s(o)=s_0$\,.
Since $\Hol(\bbE)$ is connected (because $M$ is simply connected),
$s_0$ is a fixed point with respect to the action of $\Hol(\bbE)$ if and only if
\begin{equation}
\forall A\in\hol(\bbE):A\,s_0=0\;.
    \label{eq:existence_1}
\end{equation}
The lemma follows from~\eqref{eq:hol_of_bbE} and~\eqref{eq:existence_1}.
\end{proof}
We will now apply Lemma~\ref{le:R_E_is_parallel} to deduce the following
result:

\bigskip
\begin{proposition}\label{p:properties_2}
Let a parallel isometric immersion $f:M\to N$ be given. 
For each $p\in M$
there exists a unique involutive map $I_p:\bot^1f \to\bot^1f$\,, characterized by the following properties:
% \Abstand{3}
\begin{enumerate}
\item $I_p$ is a fibrewise isometric vector bundle homomorphism along
  $\sigma^M_p$\,, i.e.\
the following diagram is commutative,
$$\begin{CD}
\bot^1f @>I_p>> \bot^1f\\
@VVV @VVV \\
M @>\sigma^M_p>> M\;
\end{CD}
$$
and for each $q\in M$ the map
$I_p|\bot^1_qf:\bot^1_qf\to\bot^1_{\sigma^M_p(q)}f$ is a linear isometry.
\item $I_p$ is parallel.
\item $I_p|\bot^1_pf$ is the identity on $\bot^1_pf$\,.
\end{enumerate}
Moreover:
\begin{enumerate}\addtocounter{enumi}{3}
\item
For every smooth geodesic line $\gamma:[-1,1]\to M$ with $\gamma(0)=p$ we have
\begin{equation}\label{eq:I_p}
I_p|\bot^1_{\gamma(-1)}f=\qmq{r.h.s.\ of~\eqref{eq:Sigma_3}}\;.
\end{equation}
\item We also have for each $q\in M$\,:
\begin{equation}
\label{eq:D_sigma_is_parallel_1}
\forall x,y\in T_qM:I_ph(x,y)=h(T\sigma^M_{p}\,x,T\sigma^M_{p}\,y)\;.
\end{equation}
\end{enumerate}
\end{proposition}
\begin{proof}
% In the following we again assume (for simplicity) 
% that $M\subset N$ is a parallel submanifold and $f=\iota^M$\,.
The uniqueness of $I_p$ follows immediately from its parallelity together with~(c).
For its existence let, we consider the origin $o=p$ and put $\sigma:=\sigma^M_o$\,. To prove
the existence of $I_o$\,, we will apply Lemma~\ref{le:R_E_is_parallel} with
$\bbE:=\rmL(\bot^1f,\sigma^*\bot^1f)$\, (where $\sigma^*\bot^1f$ is
the pullback bundle, whose fiber at $q\in M$ is given by $\bot^1_{\sigma(q)}f$).
Thus $\bbE$ is a vector bundle over $M$ with fibers
$\bbE_q=\rmL(\bot_qf,\bot_{\sigma(q)}f)$\,,
whose sections correspond in a natural way to the vector bundle homomorphisms
of $\bot^1 f$ along $\sigma$\,. To get a connection on $\bbE$\,, note that
$\nabla^\bot$ defines a connection on $\bot^1f$ (since $\bot^1f\subset\bot f$ is
a parallel vector subbundle). The pullback of this connection via
$\sigma^*$ gives a connection on
$\sigma^*\bot^1f$\,; thus we obtain the induced connection on $\bbE$\,, such that parallel
sections of $\bbE$ correspond to parallel vector bundle homomorphisms. Its
parallel displacement of an element $\ell\in\bbE_q$ along a curve $c:[0,1]\to
N$ with $c(0)=q$ is given by
\begin{equation}\label{eq:pardis_in_bbE}
\ghdisp{0}{1}{c}{\bbE}(\ell)\,\xi=\ghdisp{0}{1}{\sigma \circ c}{\bot}\circ
\ell\circ\ghdisp{1}{0}{c}{\bot}\,\xi\qmq{for all}\xi\in\bot_{c(1)}^1f\;.
\end{equation}
According to Part~(c) of Proposition~\ref{p:properties}, the curvature tensor of $\bot^1f$ (which
is the restriction of $R^\bot$ to $\bot^1f$) is a parallel tensor; thus the
curvature tensor of $\sigma^*\bot^1f$ is given by
$$R^{\bot}(T_q\sigma\,x,T_q\sigma\,y)\,\xi\qmq{for all}x,y\in
T_qM\qmq{and}\xi\in \bot_{\sigma(q)}^1f\;$$ 
is also parallel (since $\sigma$
is an isometry of $M$). Therefore, the induced curvature tensor of $\bbE$ given
for all $x,y\in T_qM,\,\ell\in \bbE_q$ by
$$R^\bbE(x,y)\,\ell=R^{\bot}(T_q\sigma\,x,T_q\sigma\,y)\circ \ell -\ell\circ R^\bot(x,y)
$$
is parallel, too.
As $\bbE_o=\rmL(\bot^1_of,\bot^1_of)$\,, we obtain for $s_0:=\Id:=\Id_{\bot^1_of}$
$$\forall x,y\in T_oM:\,R^\bbE(x,y)\,s_0=[R^{\bot}(x,y),\Id]=0\;,$$
hence Equation~\eqref{eq:existence_2} holds.
Thus there exists a unique parallel section $s$ of
$\rmL(\bot^1f,\sigma^{*}\bot^1f)$ with $s(o)=\Id$\,.
Let $I_o$ denote the corresponding vector bundle homomorphism.
To verify~\eqref{eq:I_p}, notice that $\sigma\circ\gamma$ is the inverse curve $\gamma^{-1}:t\mapsto\gamma(-t)$\,.
Because of
$$I_o(\gamma(-1))=\ghdisp{0}{1}{\gamma^{-1}}{\bbE}(\Id)\stackrel{\eqref{eq:pardis_in_bbE}}
=\ghdisp{0}{1}{\underbrace{\sigma\circ\gamma^{-1}}_{=\gamma}}{\bot f}\circ
\Id\circ\ghdisp{0}{1}{\gamma^{-1}}{\bot f}=\ghdisp{1}{-1}{\gamma}{\bot f}\;,$$
Equation~\eqref{eq:I_p} follows; we also see from the last Equation that
$I_o(p)$ is an isometry for each $p\in M$\,, since $p$ can be joint with
$o$ through some geodesic.
We have $I_o^2=\Id$ on $\bot^1_of$\,,
thus $I_o^2=\Id$ follows from the parallelity of $I_o$\,.

To prove
\eqref{eq:D_sigma_is_parallel_1}, let $c:[0,1]\to M$ be a curve
with $c(0)=p$ and $c(1)=q$\,, put $\sigma:=\sigma^M_p$\,.
Let $X,Y$ be parallel sections of $TM$ along $c$ with $X(1)=x$
and $Y(1)=y$\,. Consider the two sections $S_1$ and $S_2$ of $\bot^1f$ along the curve
$\sigma\circ c$ defined by $S_1(t):=I_{p}(h(X(t),Y(t))$ and $S_2(t):=h(T\sigma\,X(t),T\sigma\,Y(t))$\,.
Using (c) and the parallelity of $h$\,, we see that
$S_1$ is a parallel section. $S_2$ is
parallel, too, because $\sigma$ is an isometry of $M$\,. Furthermore
$S_1(0)=S_2(0)$ holds, since we have (with $\tilde x:=X(0)$\,, $\tilde y:=Y(0)$):
$$
S_1(0)=I_{p}h(\tilde x,\tilde y))\stackrel{(c)}=h(\tilde x,\tilde y)=h(-\tilde
x,-\tilde y)=h(T_p\sigma\,\tilde x,T_p\sigma\,\tilde y)=S_2(0)\;;
$$
therefore $S_1=S_2$\,, in particular~\eqref{eq:D_sigma_is_parallel_1} holds.
\end{proof}

\bigskip
\begin{remark}
Even if $f$ is not 
parallel, then nevertheless it may happen that the involution $I_p$ described
above exists. But one can easily show that~\eqref{eq:D_sigma_is_parallel_1} in addition implies the parallelity
of $f$\,.
\end{remark}

\bigskip
\begin{proof}[Proof of Theorem~\ref{th:wes}]
The uniqueness of the 
described map on $\osc f$ follows immediately from Properties~(b) and~(c) described in Theorem~\ref{th:wes}.
To prove its existence,
we consider for each $p\in M$ the unique vector bundle isomorphism of $\osc f$
given by
\begin{align}\label{eq:Sigma}
\forall q\in M,\forall x\in T_qM,\xi\in\bot^1_qf:\Sigma_p(x+\xi)=T_p\sigma^M_p(x)+I_p(\xi)\;,
\end{align}
where $I_p$ was defined in Proposition~\ref{p:properties_2}.
Then $\Sigma_p$ is an involution of $\osc f$ and a fibrewise isometric, split-parallel
(cf. Definition~\ref{de:split-parallel}) vector bundle homomorphism  along
$\sigma_p^M$ according to~\eqref{eq:Sigma} and Parts~(a),~(b) of
Proposition~\ref{p:properties_2}. Furthermore,~\eqref{eq:Sigma} combined with
Part~(c) of Proposition~\ref{p:properties_2} and the equality
$T_p\sigma_p^M=-\id_{T_pM}$ implies Part~(c) of Theorem~\ref{th:wes},
whereas~\eqref{eq:Sigma_restricted_to_TM} and~\eqref{eq:Sigma_2} follow from~\eqref{eq:Sigma}
combined with the well known facts that we have
$\sigma_p^M(\gamma(-1))=\gamma(1)$ and
$$
T_{\gamma(-1)}\sigma_p^M=-\ghdisp{-1}{1}{\gamma}{M}\;.
$$
\eqref{eq:Sigma_is_split-parallel}
is an immediate consequence of~\eqref{eq:Sigma} combined with~\eqref{eq:D_sigma_is_parallel_1}.
\eqref{eq:Sigma_3} follows immediately from~\eqref{eq:Sigma} in combination
with Part~(d) of Proposition~\ref{p:properties_2}. 
It remains to establish Assertion~(b) of Theorem~\ref{th:wes}:

For this, Equation~\eqref{eq:first_Gauss} (the Gau{\ss}-Weingarten equation) and the split-parallelity of $\Sigma_p$ implies that $\Sigma_p$ is
$\nabla^N$-parallel if and only if for all $q\in M$
\begin{equation}
\label{eq:D_sigma_is_parallel_3}
\forall x\in T_qM,v\in\osc_qf:\Sigma_{p}(\fetth(x)\,v)=\fetth(T\sigma^M_{p}\,x)(\Sigma_p\,v)\;.
\end{equation}
The result follows, since Equations~\eqref{eq:D_sigma_is_parallel_3}
and~\eqref{eq:Sigma_is_split-parallel} are equivalent as a consequence of
Lemma~\ref{le:splitting}.

Now suppose that $M$ is even a symmetric submanifold of $N$ with extrinsic
symmetries $\sigma_p^\bot$ ($p\in M$), according to Definition~\ref{de:extrinsically_symmetric}. Then we have 
for all $x,y\in T_p M:\,T\sigma_p^\bot (h(x,y))=h(T\sigma_p^\bot(x),T\sigma_p^\bot (y))$\,, hence $T\sigma_p^\bot
(\osc M)\subset \osc M$
and thus $\tilde\Sigma_p:=T\sigma_p^\bot|\osc M:\osc M\to \osc M$ satisfies Properties (a)-(c) stated in
Theorem~\ref{th:wes}. Hence $\Sigma_p=\tilde\Sigma_p$\,, by uniqueness.
\end{proof}

\subsection{Homogeneity of $\bot^1f$}
{\bf In this section, $\mathbf{M}$ is a simply connected Riemannian symmetric
  space which is isometric to a Riemannian product $M_1\times\cdots\times
  M_k$ of {\em irreducible} \/symmetric spaces $\mathbf{M_i}$ and
$\mathbf{f:M\to N}$ is a parallel isometric immersion}\,.  
We aim to prove that $\bot^1f$ is a homogeneous vector
bundle over $M$\,. Let $\fraki(M)=\frakk^M\oplus\frakp^M$ denote the Cartan
decomposition, and let $\Sym(M)$ denote the subgroup of $\Iso(M)$ generated by its geodesic
symmetries $\sigma_p^M$\,, where $p$ ranges over $M$\,.
One can show that $\Sym(M)$ is actually a Lie subgroup of $\Iso(M)$ with
$\Iso(M)^0\subset\Sym(M)$ (in case $M$ is irreducible, this fact is explained in~Sec.~3.3. of~\cite{J1}).

\bigskip
\begin{definition}\label{de:homogeneous}
We will call a vector bundle $\bbE$ over $M$ a
\emph{homogeneous vector bundle} if there exists an action
$\alpha:\Iso(M)^0\times\bbE\to\bbE$ by vector bundle isomorphisms such that the
bundle projection of $\bbE$ is equivariant.
\end{definition}

In the above situation, we consider $\osc f$ as a vector bundle over $M$
equipped with the connection $\nabla^{\osc f}$ described at the
beginning of Section~\ref{sec:symmetric_strip}.

\bigskip
\begin{proposition}\label{p:homogeneity_of_the_first_normal_bundle}

\begin{enumerate}
\item There exists a natural action $\alpha:\Sym(M)\times\osc f\to\osc f$ where
  $\Sym(M)$ acts through isometric, parallel vector bundle isomorphisms, characterized as follows: For each 
point $p$ of $M$ we have \begin{equation}\label{eq:alpha_1}
\forall v\in\osc_pf:\alpha(\sigma_p^M,v)=\Sigma_p(v)\;.
\end{equation}
\item $\alpha$ splits into two actions on $TM$ resp.\ on $\bot^1f$ (denoted by $\alpha^\top$ resp.\ by $\alpha^\bot$), i.e.
\begin{equation}\label{eq:alpha_2}
\forall g\in\Sym(M),x+\xi\in T_{p}M\oplus\bot^1_{p}f:\;\alpha(g,x+\xi)=\alpha^\top(g,x)+\alpha^\bot(g,\xi)\in TM\oplus\bot^1f\;;
\end{equation} 

and we have for all $g\in\Sym(M), x\in T_pM$
\begin{equation}\label{eq:alpha_top}
\alpha^\top(g,x)=T_pg\,x\;.
\end{equation}
Furthermore, the second fundamental form
$h$ is $\alpha$-invariant in the following sense: 
\begin{equation}\label{eq:h_is_invariant}
\forall g\in\Sym(M), x,y\in T_pM: \alpha^\bot(g,h(x,y))=h(\alpha(g,x),\alpha(g,y))\;.
\end{equation}
\item $\bot^1f$ is a homogeneous vector bundle over $M$ via the action of
  $\alpha^\bot$ restricted to $\Iso(M)^0$\,. 
\item
One can also show that the normal connection on $\bot^1f$ is 
the canonical connection induced by the Cartan decomposition 
as described in Section 2.1 of~\cite{J1} (without proof). 
\end{enumerate}
\end{proposition}
\begin{proof}
Put $G:=\Sym(M)$\,, and let $\tilde G$ denote the subgroup of vector bundle
isomorphisms on $\osc f$ generated by
all $\Sigma_p$ with $p\in M$ (see~Theorem~\ref{th:wes}); thus we have the natural action $\tilde\alpha:\tilde
G\times \osc f\to\osc f$ and a surjective group homomorphism
$\pi:\tilde G\to G$ such that 
$\pi(\Sigma_p)=\sigma_p^M$ for each $p\in M$\,, 
hence by Equation~\eqref{eq:Sigma_restricted_to_TM}
\begin{equation}\label{eq:tilde_alpha_g_restricted_to_TM}
\forall g\in\tilde G:\tilde\alpha_g|TM=T\pi(g)\;.
\end{equation}
Moreover, by means of~\eqref{eq:Sigma_is_split-parallel} we have for arbitrary $g\in \tilde G$ 
\begin{equation}\label{eq:h_is_equivariant}
\forall x,y\in TM: \tilde\alpha(g,h(x,y))=h(T\pi(g)\,x,T\pi(g)\,y)\;.
\end{equation}
which implies that $\pi$ is also injective. Therefore $\pi$ is
an isomorphism; thus we may pointwise define $\alpha$ 
via $\alpha_{\pi(g)}=\tilde\alpha_g$ for $g\in\tilde G$\,. 
\eqref{eq:tilde_alpha_g_restricted_to_TM} and~\eqref{eq:h_is_equivariant} imply~\eqref{eq:alpha_2}-\eqref{eq:h_is_invariant}.
It remains to show that $\alpha$ is differentiable. By means of~\eqref{eq:alpha_2}
and~\eqref{eq:alpha_top}, this will be clear if $\alpha^\bot$ is
differentiable. For this,
let an arbitrary (differentiable) section $\xi$ of $\bot^1f$ be
given. For each point $p\in M$ there exists an open neighbourhood $U$ of $p$ in
$M$\,, vector fields $X_1,\ldots,X_k,Y_1,\ldots,Y_k$ on $U$ and $C^\infty$-functions $\lambda_1,\ldots,\lambda_k$ on $U$ such
that 
$\xi|U=\sum_{i=1}^k\lambda_ih(X_i,Y_i)$\,. Then we have
$$
\alpha^\bot\big(g,\xi\big)|U=\sum_{i=1}^k\lambda_i\,h(Tg\,X_i,Tg\,Y_i),
$$
which is a differentiable function on $G\times U$\,. It follows that
$\alpha^\bot$ is differentiable.

For (c): Since $G$ is generated by the reflections $\sigma_p^M$\,, the
equivariance of $\alpha^\bot$ follows from~\eqref{eq:Phi} combined with the
construction of $\alpha$\,. Thus $\bot^1f$ is a homogeneous vector bundle over
$M$\,. 
\end{proof}
It is planed to investigate parallel isometric immersions $f:M\to N$ defined on a
symmetric space $M$ as above in a forthcoming paper~\cite{J2}.
\section{Proof of Theorem~\ref{th:hol}}
\label{s:hol}
In this section, we will give the proof of Theorem~\ref{th:hol}.
Let $f:M\to N$ be a parallel isometric immersion defined on a simply connected
symmetric space $M$\,. Without loss of generality, we can assume that $M$ is simply connected,
for the following reason: Let $\tau:\hat M\to M$ denote the universal covering, and
consider the isometric immersion $\tilde f:=f\circ\tau$ and the corresponding
holonomy group $\Hol(\tilde f^*TN)$ with respect to some point $\hat
o\in\tau^{-1}(o)$\,. Then it is well known that the connected components of $\Hol(\tilde f^*TN)$ and $\Hol(f^*TN)$
are equal, and thus the holonomy Lie algebras $\hol(\tilde f^*TN)$ and
$\hol(f^*TN)$ are equal, too\;. Moreover, for the sake of an easier notation, we
assume that $M$ is a submanifold of $N$ and $f=\iota^M$\,.

\bigskip
\begin{proof}[For Part~(c)]
Let us first lead the discussion on the level of the corresponding Holonomy groups. 
According to~\eqref{eq:Hol(TN|M)},~\eqref{eq:Hol(osc_M)} combined with~\eqref{eq:parisp_osc},
by $\Hol(f^*TN)\to \Hol(\osc f), g\mapsto g^{\osc f}$ is
defined a surjective Lie group homomorphism. This map is even an isomorphism, which is seen
as follows: Suppose $g^{\osc f}=\Id$ on $\osc_of$ for some $g\in\Hol(f^*TN)$\,.
Using the $\nabla^N$ parallelity of $R^N$\,, we have
$$
\forall u,v,w\in T_oN:g(R^N(u,v)\,w)=R^N(g\,u,g\,v)(g\,w)\;;
$$
hence the linear space $V:=\Menge{v\in T_oN}{g\,v=v}$ is curvature
invariant. Therefore, because $\osc_of\subset V$ by assumption, $f$ maps into the totally
geodesic submanifold defined by $V$\,, according to Theorem~\ref{th:Dombi};
thus $g=\Id$ by the fullness of $f$\,. 
Switching to the level of the Lie algebras, the result hence follows.
\end{proof}

\bigskip
\begin{proof}[For Part~(a)]
Let  $\sigma^\bot:\osc_of\to\osc_of$ denote the linear reflection in $\bot^1_of$\,.
We have to show that 
\begin{equation}\label{eq:hol_2}
\Ad(\sigma^\bot)(\hol(\osc f))=\hol(\osc f)\;.
\end{equation}
Let $\Sigma_o$ denote the symmetry of $\osc f$ at the point
$o$ described in Theorem~\ref{th:wes}, and let $c:[0,1]\to M$ be a loop with $c(0)=o$\,.
Remember that $\Sigma_o$ is a $\nabla^{\osc f}$-parallel vector bundle
isomorphism of $\osc f$ along $\sigma_o^M$ with $\Sigma_o|\osc_of=\sigma^\bot$\,, in accordance with
Theorem~\ref{th:wes}; hence
\begin{equation}
\label{eq:D_sigma_o_circ_F}
\sigma^\bot\circ\ghdisp{0}{1}{c}{\osc f}=\ghdisp{0}{1}{\sigma^M_o\circ c}{\osc
f}\circ \sigma^\bot\;.
\end{equation}
From the last line we
conclude that $\Hol(\osc f)$ is invariant by group conjugation with
$\sigma^\bot$\,; thus~\eqref{eq:hol_2} holds.
\end{proof} 

\bigskip
\begin{proof}[For Equation~\eqref{eq:[fetth(x),_]}]
Let $\ad:\so(\osc_of)\to\End\big(\so(\osc_of)\big),
A\mapsto[A,\,\cdot\,]$ denote the adjoint representation of
$\so(\osc_of)$\,; thus~\eqref{eq:[fetth(x),_]} is equivalent to
\begin{equation}\label{eq:ad(fetth(x))}
\forall x\in T_oM:\ad\big(\fetth(x)\big)\big(\hol(\osc f)\big)\subset\hol(\osc f)\;.
\end{equation}
Let $x\in T_oM$ and $\gamma$ be the geodesic of $M$ with $\gamma(0)=o,\;\dot \gamma(0)=x$\,,
and let $\Theta_\gamma(t)\;(t\in \R)$ denote the family of vector bundle
isomorphisms on $\osc f$ along $\theta_\gamma(t)$ from Definition~\ref{de:Phi}.
Because $\Theta_\gamma(t)$ is a $\nabla^{\osc f}$-parallel vector bundle
isomorphism of $\osc f$ along $\theta(t)$ for each $t\in\R$ (Corollary~\ref{co:Phi}), we obtain for each loop $c:[0,1]\to M$ with $c(0)=o$
\begin{equation}\label{eq:Phi_c}
\ghdisp{0}{t}{\gamma}{sp}\circ\ghdisp{0}{1}{c}{\osc f}\circ\ghdisp{t}{0}{\gamma}{sp}|\osc_{\gamma(t)}f
\stackrel{\eqref{eq:Phi_gleich_phi}}{=}\Theta_\gamma(t)\circ\ghdisp{0}{1}{c}{\osc
  f}\circ \Theta_\gamma(t)^{-1}|\osc_{\gamma(t)}f
=\ghdisp{0}{1}{\theta_\gamma(t)\circ c}{\osc f}\;;
\end{equation}
and therefore (by virtue of~\eqref{eq:osc_is_parallel_1},~\eqref{Eq:Move2} and~\eqref{eq:mu_c_2})
\begin{equation}\label{eq:mu_c_3}
\mu_\gamma(t)\circ\ghdisp{0}{1}{c}{\osc f}\circ
\mu_\gamma(t)^{-1}=\ghdisp{0}{1}{c_t}{\osc f}\;,
\end{equation}
where $c_t$ means the loop at
$o$  defined by going first along $\gamma$ from $o$ to $\gamma(t)$\,,
then along the loop $\theta_\gamma(t)\circ c$ centered at $\gamma(t)$ and then back from $\gamma(t)$ to
$o$ along the inverse curve $\gamma^{-1}$\,.
In accordance with Equations~\eqref{eq:Hol(TN|M)} and~\eqref{eq:Hol(osc_M)},
this shows that l.h.s.\ of~\eqref{eq:mu_c_3} is an
element of $\Hol(\osc f)$\,. Therefore,
we have for each $t\in \R$
\begin{equation}\label{eq:phi_acts_on_hol(osc_M)}
\mu_\gamma(t)\circ \Hol(\osc f)\circ
\mu_\gamma(t)^{-1}\subset\Hol(\osc f)\;,\qmq{hence}\Ad(\mu_\gamma(t))(\hol(\osc f))\subset\hol(\osc f)
\end{equation}
(where $\Ad$ means the adjoint representation of $\SO(\osc_of)$).
Since $\mu_\gamma(t)\stackrel{\eqref{eq:mu_1}}{=}\exp(t\,\fetth(x))$\,,
Equation~\eqref{eq:ad(fetth(x))} follows by taking the derivative
in~\eqref{eq:phi_acts_on_hol(osc_M)} with respect to $t$ at $t=0$\,.
\end{proof} 

\bigskip
\begin{proof}[For Part~(b)] 
In the following, the simple relations between the parallel displacement resp.\ the
curvature tensor of $\osc f$ and $f^*TN$ described in~\eqref{eq:osc_is_parallel_1}-\eqref{eq:R_osc_1} will be used
without further reference. Taking into account~\eqref{eq:fetth_in_osc} and~\eqref{eq:fundamental_2}, we can define
the following linear subspaces of $\so(\osc_of)$\,:
\begin{align*}
&\frakj_i:=\Spann{[\fetth(x_{i}),[\fetth(x_{i-1}),\ldots,[\fetth(x_1),R^{\osc f}(y_1,y_2)],\ldots,]}{x_1,\ldots,x_i\in
  T_oM,y_1,y_2\in T_oM}\;,\\
&\frakj:=\sum_{i=0}^3\frakj_i\;,\\
&\frakj_+:=\text{r.h.s.\ of}~\eqref{eq:hol_plus}\;,\\
&\frakj_-:=\text{r.h.s.\ of}~\eqref{eq:hol_minus}\;.
\end{align*}
 Because of Equations~\eqref{eq:R_osc_2} and~\eqref{eq:fundamental_1}, we have $\frakj_+\subset\so(\osc_of)_+$\,; and hence
 $\frakj_-\subset\so(\osc_of)_-$\,, according to~\eqref{eq:fetth_in_osc} and the
 rules for $\Z/2\Z$ graded algebras.
Let us now see that $\hol(\osc f)=\sum_{i=0}^3\frakj_i=\frakj_+\oplus\frakj_-$ holds; the proof will be divided into three steps.

First step.
Let us see that we have
$\frakj_+\oplus\frakj_-\subset\frakj\subset\hol(\osc f)$\,:
As a consequence of the Theorem of Ambrose/Singer, we have
$\frakj_0\subset\hol(\osc f)$ and hence $\frakj_i\subset\hol(\osc f)$ for each
$i=0,\ldots,3$\,, according to~\eqref{eq:[fetth(x),_]}; thus $\frakj\subset\hol(\osc f)$\,.
Moreover, Proposition~\ref{p:fundamental} implies that
$$
\forall\xi,\eta\in
\bot^1_of:\big(R^N(\xi,\eta)\big)^{\osc}\in\frakj_0+\frakj_2\;;
$$
thus also
$
\frakj_+\subset\frakj_0+\frakj_2$ and $\frakj_-\subset\frakj_1+\frakj_3\;.
$

Second step.
I claim that $\frakj_+\oplus\frakj_-\subset\so(\osc_of)$ is a vector space
invariant by
$\ad(\fetth(x))$ for each $x\in T_oM$\,:
It suffices to show that $$
[\fetth(x),\frakj_-]\subset \frakj_+\;,
$$
which means for all $z_1,z_2\in T_oM$\,, $\xi,\eta\in\bot^1_of$\,:
\begin{align}\label{eq:tildex}
&[\fetth(x),[\fetth(y),\big(R^N(z_1,z_2)\big)^{\osc}]]\in\frakj_+\qmq{and}\\
&\label{eq:fundamental_3}
[\fetth(x),[\fetth(y),(R^N(\xi,\eta))^{\osc}]]\in\frakj_+\;.
\end{align}
\eqref{eq:tildex} holds because of Proposition~\ref{p:fundamental}. For~\eqref{eq:fundamental_3} choose $v\in \osc_of$\,; then
by means of Lemma~\ref{le:more_general} (applied to~\eqref{eq:parallel_2}) we get
\begin{align}\notag
&[\fetth(x),[\fetth(y),R^N(\xi,\eta)]]\,v
=R^N\underbrace{(\fetth(x)\,\fetth(y)\,\xi,\eta)}_{\in U\times U}\,v
+R^N\underbrace{(\xi, \fetth(x)\,\fetth(y)\,\eta)}_{\in U\times U})\,v\\
\label{eq:fetth_4}
&+R^N\underbrace{(\fetth(x)\,\xi,\fetth(y)\,\eta)}_{\in W\times W}\,v
+R^N\underbrace{(\fetth(y)\,\xi,\fetth(x)\,\eta)}_{\in W\times
  W}\,v\;,
\end{align}
with $U:=\bot^1_of$\,, $W:=T_oM$\,; which proves~\eqref{eq:fundamental_3}.

Third step.
$\hol(\osc f)\subset \frakj_+\oplus\frakj_-$ is finally proved as follows: 
By virtue of the Theorem of Ambrose/Singer, the vector space
$\hol(\osc f)$ is generated by elements of the form
\begin{equation}\label{eq:Ambrose}
\ghdisp{1}{0}{c}{\osc f}\circ R^{\osc f}(y_1,y_2)\circ
\ghdisp{0}{1}{c}{\osc f}
\end{equation}
for various curves $c:[0,1]\to M$ with $c(0)=o$ and $y_1,y_2\in T_{c(1)}M$\,.
Therefore, given such a curve $c$ and $y_1,y_2\in
T_{c(1)}M$\,, we introduce 
$$\forall t\in[0,1]:\tilde R(t):=\ghdisp{t}{0}{c}{\osc f}\circ R^{\osc f}(y_1,y_2)\circ
\ghdisp{0}{t}{c}{\osc f}\;\in\so(\osc_of)\;.$$
Of course, by the previous, it suffices to show that
\begin{equation}\label{eq:R_2}
\forall t\in[0,1]:\tilde R(t)\in \frakj_+\oplus\frakj_-\;.
\end{equation}

For this: Let $\mu_c$ be the function defined in Equation~\eqref{Eq:Move2}.
From Lemma~\ref{le:R(t)} we get
\begin{align}
&\tilde R(t)=\Ad(\mu_c(t))\,R^{\osc f}(\tilde y_1,\tilde y_2)\;,
    \label{eq:hol_3}
\end{align}
with $\tilde y_i:=\ghdisp{t}{0}{c}{sp}\,y_i$ for $i=1,2$\,. Introduce the linear space
\begin{equation}\label{eq:frakn}
\frakn_\ad(\frakj_+\oplus\frakj_-):=\Menge{A\in\so(\osc_of)}{\ad(A)(\frakj_+\oplus\frakj_-)\subset\frakj_+\oplus\frakj_-}\;,
\end{equation}
which is actually the Lie algebra of the Lie subgroup of
$\SO(\osc_of)$ given by
\begin{equation}\label{eq:rmN}
\rmN_\Ad(\frakj_+\oplus\frakj_-):=\Menge{g\in\SO(\osc_of)}{\Ad(g)(\frakj_+\oplus\frakj_-)=\frakj_+\oplus\frakj_-}\;.
\end{equation}
By means of the second step, we have
$\fetth(X(t))\in\frakn_\ad(\frakj_+\oplus\frakj_-)$ for
each $t\in \R$ (where $X:[0,1]\to T_oM$
denotes the function defined by~\eqref{eq:X(t)}), and therefore the left
invariant vector field $\tilde X$ defined on the Lie group $\SO(\osc_of)$ by $\forall g\in\SO(\osc_of):\tilde X_t(g):=g\circ \fetth(X(t))$  is tangential to the submanifold
$\rmN_\Ad(\frakj_+\oplus\frakj_-)$\,. By means of~\eqref{Eq:Move3}, the curve
$\mu_c$ solves the ODE
$$
\dot\mu_c(t)=\tilde X_t(\mu_c(t))\qmq{with}\mu_c(0)=\Id\;.
$$
Thus we find that in fact $\mu_c$ is a curve in $\rmN_\Ad(\frakj_+\oplus\frakj_-)$\,.
Since $\tilde y_i\in T_oM$ for $i=1,2$ (in accordance with Definition~\ref{de:split-parallel}), we
moreover have
$R^{\osc f}(\tilde y_1,\tilde y_2)\in\frakj_+$\,. From the
previous, we finally conclude that
$\forall
t\in[0,1]:t\mapsto\Ad(\mu_c(t))\,R^{\osc f}(\tilde
y_1,\tilde y_2)$
actually describes a curve into $\frakj_+\oplus\frakj_-$\,, which together
with~\eqref{eq:hol_3} proves~\eqref{eq:R_2}.
\end{proof}

\section{The extrinsic holonomy Lie algebra of a full symmetric submanifold\ldots}
\label{s:extrinsically_hol}
In this section, $M^m$ is a full symmetric submanifold 
of a simply connected, irreducible symmetric space
$N$\,, and $o\in M$ is some origin. We will now
calculate as explicitly as possible the extrinsic holonomy Lie algebra of
$M$\,.

\bigskip
\begin{proposition}\label{p:frakk_gleich_hol}
In the above situation, the Lie algebras $\frakk$ and $\hol(N)$ are isomorphic via $\pi_2$\,.
\end{proposition}

\begin{proof}
It is well known that $\pi_2$ is a faithful representation of $\frakk$ on $T_oN$\,;
thus it is sufficient to verify that $\pi_2(\frakk)=\hol(N)$ holds.
By the Theorem of Ambrose/Singer we have
$\hol(N)=\Spann{R^N(x,y)}{x,y\in T_oN}$\,.
Since $N$ is irreducible, the Lie algebra $\fraki(N)$ is semisimple
(cf.~\cite{He},~Ch.\,V,~Prop.\,4.2) and hence we can apply~\cite{He},~Ch.\,V,~Part~(iii) of~Theorem~4.1\,.
\end{proof}

In accordance with Definition~\ref{de:extrinsically_symmetric}, 
let $\sigma^\bot_o\in\Iso(N)$ be the corresponding extrinsic symmetry of $M$ at $o$\,.
Since $\sigma_o^\bot$ is an isometry of $N$ with
$\sigma_o^\bot(o)=o$ and $T_o\sigma^\bot_o=\sigma^\bot$\,, we have
\begin{equation}\label{eq:symmetric_pair}
\sigma^\bot\circ\Hol(N)\circ\sigma^\bot=\Hol(N)\qmq{and}\Ad(\sigma^\bot)\,\hol(N)=\hol(N)\;.
\end{equation}
(where $\sigma^\bot:T_oN\to T_oN$ denotes the linear reflection in
$\bot_oM$); therefore the splitting $T_oN=T_oM\oplus\bot_oM$ induces the splitting
\begin{align}
\label{eq:decomposition_of_hol(N)}
&\hol(N)=\hol(N)_+\oplus\hol(N)_-\qmq{with}\\
&\hol(N)_\pm:=\hol(N)\cap\so(T_oN)_\pm\;.\label{eq:hol(N)_pm}
\end{align}
\bigskip
\begin{lemma} We have
\begin{align}\label{eq:hol_von_N_plus}
&\hol(N)_+=\Spann{R^N(x,y)}{x,y\in
T_oM}+\Spann{R^N(\xi,\eta)}{\xi,\eta\in\bot_oM}\qmq{and}\\
\label{eq:hol(N)_-}
&\hol(N)_-=\Spann{R^N(x,\xi)}{x\in T_oM,\xi\in\bot_oM}\;.
\end{align}
\end{lemma}
\begin{proof}
On the one hand, $\sigma^\bot R^N(u,v)\sigma^\bot=R^N(\sigma^\bot\,u,\sigma^\bot\,v)$ for all $u,v\in T_oN$\,,
thus $R^N(u,v)\in\so(T_oN)_+$ (resp.\ $R^N(u,v)\in\so(T_oN)_-$)
if $u$ and $v$ are both contained in
$T_oM$ or both  in $\bot_oM$ (resp.\ if $u\in T_oM$ and $v\in\bot_oM$).
On the other hand, $\hol(N)=\Spann{R^N(x,y)}{x,y\in T_oN}$ by the
Theorem of Ambrose/Singer (since $R^N$ is a parallel tensor).
\end{proof}

\bigskip
\begin{proposition}\label{p:frakk_bot}
We have $\osc M=TN|M$\,. Consequently, we can introduce the splitting
$\hol(TN|M)=\hol(TN|M)_+\oplus\hol(TN|M)_-$\,, in accordance
with~\eqref{eq:splitting_of_hol}. The Lie algebra $\hol(TN|M)_+$ coincides with $\hol(N)_+$\,.
Moreover, $\hol(TN|M)_-$ is an  $\ad(\hol(N)_+)$-invariant linear subspace of
$\hol(N)_-$\, i.e.\ we have \begin{equation}\label{eq:Invariance_6.0}
[\hol(N)_+,\hol(TN|M)_-]\subset \hol(TN|M)_-\;.
\end{equation}
\end{proposition}
\begin{proof}
$M$ is a 1-full, parallel submanifold of $N$ by virtue of
Theorem~\ref{th:result_1}, thus $\osc M=TN|M$ holds. Comparing~\eqref{eq:hol_plus} and~\eqref{eq:hol_von_N_plus},
we hence see that we have $\hol(TN|M)_+=\hol(N)_+$\,. Since $\hol(TN|M)$ is a Lie algebra, the last assertion now follows
from the rules of $\Z/2\Z$-graded Lie algebras.
\end{proof}

Let $\Hol(N)_+:=\Menge{g\in\Hol(N)}{\sigma^\bot\circ g\circ\sigma^\bot=g}$\;;
then, according to~\eqref{eq:symmetric_pair}, $\big(\Hol(N),\Hol(N)_+\big)$ is
a symmetric pair (in the sense of~\cite{He},~Ch.\,4,~\S~3),
and, in accordance with~\eqref{eq:hol(N)_pm}, the Lie algebra of $\Hol(N)_+$ is given by $\hol(N)_+$\,.
Moreover, with respect to the natural action of $\SO(T_oN)$ on the
Grassmannian $\rmG_m(T_oN)$\,, the isotropy subgroup of $\Hol(N)$ at $T_oM$ is
given by $\Hol(N)_+$\,; hence the quotient space
$L:=\Hol(N)/\Hol(N)_+$ is equipped with a natural inclusion $L\hookrightarrow
\rmG_m(T_oN)\,,\; g\circ\Hol(N)_+\mapsto T_og(T_oM)$\,.
This maps $L$ onto a totally geodesic submanifold of the
symmetric space $\rmG_m(T_oN)$\,; then the metric on the tangent space
$T_{[e]}\,L\cong\hol(N)_-$ is given by $\g{A}{B}=-1/2\cdot\trace(A\circ B)$\,.
In this way, $L$ becomes a
Riemannian symmetric space, and, moreover,~\eqref{eq:decomposition_of_hol(N)} defines an orthogonal
symmetric Lie algebra (in the sense of~\cite{He},~Ch.\,V,~\S~1) such that the symmetric pair
$\big(\Hol(N),\Hol(N)_+\big)$ is associated therewith
(in the sense of~\cite{He},~Ch.\,IV, definition preceeding Prop.\,3.6). 

\setcounter{footnote}{0}
\bigskip
\begin{lemma}\label{le:about_hol}
% \Abstand{3}
\begin{enumerate}
\item
Let $\hat L$ denote the universal covering space of $L$\,. 
Then there exist 
symmetric spaces $L_1$ and $L_2$ such that $\hat L\cong L_1\times L_2$ with
$T_{[e]}L_1\cong\hol(TN|M)_-$ and  $T_{[e]}L_2\cong\big(\hol(TN|M)_-\big)^\bot$\,. In particular, if neither $\hol(TN|M)_-=\hol(N)_-$ nor
$\hol(TN|M)=\{0\}$\,, then $L$ is reducible\footnote{In accordance with~\cite{BCO},~p.\,290,
  we use the following convention: A Riemannian manifold $M$ is called ``reducible'' if its universal covering splits as
a (non-trivial) product of two Riemannian spaces; otherwise $M$ is called
``irreducible''.}.
\item For each subspace $V\subset\so(T_oN)$ we introduce its centralizer in
  $\so(T_oN)$\,, via $$\frakc(V):=\Menge{A\in\so(T_oN)}{\forall B\in V:A\circ
    B=B\circ A}\;.$$ If $\dim\big(\frakc(\hol(N)_+)\cap\hol(N)_-\big)<m$\,, then  $\hol(TN|M)_-\neq\{0\}$\,.
\item Suppose that $N$ is a Hermitian symmetric space and that $M\subset N$
is a Lagrangian submanifold; let $\j$ denote the complex structure of
$T_oN$\,. Then $\j$ is orthogonal to $\hol(TN|M)$\,.
\item Suppose that $N^{4n}$ is a quaternionic K{\"a}hler symmetric space with $n\geq 2$
  and that $M^{2n}$ is a totally complex submanifold. Let a canonical
  basis $\{\i,\j,\k\}$ of the quaternionic structure of $T_oN$ be given such that $\i(T_oM)=T_oM$ and $\j(T_oM)=\bot_oM$
  holds. Then both $\j$ and $\k$ are orthogonal to $\hol(TN|M)$\,.
\end{enumerate}
\end{lemma}
\begin{proof}
For~(a): As before, we consider the orthogonal symmetric Lie algebra
$\frakg:=\hol(N)$\,. Suppose first that $\hat L$ is of compact type.
Then, in accordance with Proposition~\ref{p:frakk_bot}
 $$[\hol(N)_+,\hol(TN|M)_-]\subset \hol(TN|M)_-\qmq{and}[\hol(N)_+,\big(\hol(TN|M)_-\big)^\bot]\subset\big(\hol(TN|M)_-\big)^\bot\;;$$
therefore, Proposition~\ref{p:frakk_gleich_hol} in combination with~\cite{He},~Ch.\,V,~Part~(i)
of~Theorem~4.1, (applied to $L$) shows that both $\hol(TN|M)_-$ and
$\big(\hol(TN|M)_-\big)^\bot$ are $\hol(L)$-invariant
subspaces of $\hol(N)_-\cong T_{[e]}L$\,. Thus the result follows from the
decomposition theorem of de\,Rham.

In the general case, note that  the sectional curvature of $L$ is non-negative (because $L$ is
totally geodesically embedded in $\rmG_m(T_oN)$), hence,
according to~\cite{He},~Theorem~3.1, in combination with Proposition~4.2, $\hat
L$ is a product of a Euclidian space and a symmetric space of compact type
(in fact, it may happen that $\hat L$ splits off a Euclidian factor; for
example, if $N$ is a complex space form and $M$ is a Lagrangian
submanifold of $N$ -- cf.\ the proof of Theorem~\ref{th:hol_for_symmetric_M} in the next section); however,
using a decomposition theorem for orthogonal symmetric Lie algebras
(see~\cite{He},~Theorem~1.1), and switching to the level of symmetric pairs,
we easily reduce the problem to the case when already $L$ is of compact type (cf.\ the proof of~\cite{He},~Proposition~4.1).

For~(b):
Assume that $\hol(TN|M)_-=\{0\}$ holds\,. Then we have 
$[\fetth(x),\hol(N)_+]=\{0\}$ according
to~\eqref{eq:hol_minus} and Proposition~\ref{p:frakk_bot}, hence
$\fetth(x)\in\frakc(\hol(N)_+)$ for each $x\in T_oM$\,. Furthermore, by virtue
of~\eqref{eq:Eschenburg} combined with
Proposition~\ref{p:frakk_gleich_hol}, 
for each $x\in T_oM$ we
have $$\fetth(x)\in\pi_2(\frakk)\cap\so(T_oN)_-\stackrel{\eqref{eq:hol(N)_pm}}{=}\hol(N)_-\;.$$ 
Thus  $\fetth(T_oM)$ is actually a subspace of
$\frakc(\hol(N)_+)\cap\hol(N)_-$ and therefore
$\dim\big(\frakc(\hol(N)_+)\cap\hol(N)_-\big)\geq m$\,, as a consequence of
Part~(c) of Theorem~\ref{th:result_1}.

For~(c): If $M\subset N$ is a Lagrangian submanifold, then $J_p$ maps the tangent space of $T_pM$ onto the normal space
$\bot_pM$ and vice versa; thus we have $J_p\in\so(T_pN)_-$ in accordance with~\eqref{eq:odd}, whereas the curvature invariance of $T_pM$
implies that $R^N(x,y)\in\so(T_pN)_+$ for all points $p\in M$ and $x,y\in T_pM$\,; hence
$\trace(J_p\circ R^N(x,y))=0$\,. Therefore, using the parallelity of $J$
in combination with the Theorem of Ambrose/Singer, we see 
that $\trace(\j\circ A)=0$ holds for every $A\in\hol(TN|M)$\,, i.e. $\j\in\hol(TN|M)^\bot$\,.

For~(d): 
Let $\rmQ$ denote the $\nabla^N$-parallel subbundle of $\so(TN)$ which
defines the Quaternionic K{\"a}hler structure of $N$\,.
Since $T_oM$ is a totally complex subspace of $T_oN$\,, according to Definition
2.7 of~\cite{Ts4} there exists a canonical basis $\{\i,~\j,~\k\}$ of $\scrQ_o$
with the additional property that $\i(T_oM)=T_oM$ and
$\j(T_oM)=\bot_oM$\,. Such a canonical basis is not unique; however, if $\{\tilde\i,~\tilde\j,~\tilde\k\}$ is a second
canonical basis of $\rmQ_o$ with this property, then we have $\tilde\i=\pm\i$ and there exists some $\varphi\in [0,2\pi]$
such that $\tilde\j=\cos(\varphi)\,\j-\sin(\varphi)\,\k$ and
$\tilde\k=\sin(\varphi)\j+\cos(\varphi)\k$\,. 
Furthermore, by definition of a totally complex submanifold (see~\cite{Ts4},
Definition~2.8), the pullback bundle $\rmQ|M$ is locally spanned by three
sections $I,J,K$ which satisfy the usual quaternionic relations such that additionally $I_p(T_pM)=T_pM$ and
$J_p(T_pM)=\bot_pM$ for all $p$\,. By the previous, without loss of generality we may
assume that $I_o=\i$\,, $J_o=\j$ and $K_o=\k$
holds. Since $n\geq 2$\,, $I$
is even a $\nabla^N$-parallel section of $\rmQ|M$\,, according to Lemma~2.10 of~\cite{Ts4}\,; hence the vector
subbundle $\tilde\rmQ$ which is locally spanned by $\{J,K\}$ is even a
globally well defined, $\nabla^N$-parallel subbundle
of $\rmQ|M$\,. Furthermore, $\tilde\rmQ_p\subset\so(T_pN)_-$ for each $p\in M$ (like in the
Lagrangian case); now a proof which
uses the same ideas as for Part~(c) shows that $\j$ and $\k$ both belong to
the orthogonal complement of $\hol(TN|M)$ in $\so(T_oN)$\,.  
\end{proof}
\subsection{\ldots in an ambient symmetric space of rank 1}
\label{s:rank-one}
\begin{proof}[Proof of Theorem~\ref{th:hol_for_symmetric_M} in case the rank
  of $N$ is 1]
According to Theorem~\ref{th:Naitoh}, it suffices to consider the following  possibilities.

{\bf $\mathbf N$ is a real space form of non-vanishing sectional curvature:} Here we have $\hol(N)=\so(T_oN)$ and we claim
that $\hol(TN|M)=\so(T_oN)$ holds, unless $(n,m)=(4,2)$ and
$\dim(\hol(TN|M))=4$\,:

For this, remember that $\so(T_oN)_+\cong\so(T_oM)\oplus\so(\bot_oM)$
(see~\eqref{eq:even}), and hence straightforward considerations show that
$\frakc(\so(T_oN)_+)\cap\so(T_oN)_-=\{0\}$ holds (since $m\geq 2$ by assumption).
Thus the possibility $\hol(TN|M)_-=\{0\}$ is excluded, as a consequence of Part~(b) of Lemma~\ref{le:about_hol}.  
Furthermore, the symmetric space $L$ (described in the last section) corresponds to
the real Grassmannian $\rmG_m(T_oN)$\,, which is an irreducible
symmetric space unless $(n,m)\neq(4,2)$\,. Therefore, if $(n,m)\neq(4,2)$\,, then we
have $\hol(TN|M)=\hol(N)$ pursuant to Part~(a) of
Lemma~\ref{le:about_hol}. For $(n,m)=(4,2)$\,, $\rmG_m(T_oN)$ is a 4-dimensional
reducible symmetric space; its universal covering splits into two 2-dimensional factors. From
Part~(a) of Lemma~\ref{le:about_hol}, combined with the previous, we conclude that in
this case $\dim(\hol(TN|M)_-)\in\{2,4\}$ and $\hol(TN|M)_+=\so(T_oN)_+$
holds. We thus obtain $\dim(\hol(TN|M))\in\{4,6\}$\,; moreover, in case $\dim(\hol(TN|M))=6$ we necessarily
have $\hol(TN|M)=\so(T_oN)$\,. 

\bigskip
{\bf $\mathbf N$ is a complex space form and $\mathbf M$ is a complex submanifold:} Here we have
$\hol(N)=\fraku(T_oN)=\R\,\j\oplus\su(T_oN)$\,, where $\j$ denotes the complex
structure of $T_oN$\,. We will show that always
$\hol(TN|M)=\fraku(T_oN)$ holds, as follows:

We have $\hol(N)_+=\fraku(T_oN)\cap\so(T_oN)_+\cong\su(T_oM)\oplus\su(\bot_oM)\oplus\j$
and thus we easily verify that
$\frakc(\hol(N)_+)\cap\so(T_oN)_-=\{0\}$ holds. Moreover, $L$ is isomorphic to the
Grassmannian manifold of complex $m$-planes in $T_oN$\,, which is an
irreducible symmetric space; therefore, by combining Parts~(a) and~(b) of Lemma~\ref{le:about_hol}, we obtain $\hol(TN|M)=\hol(N)$\,.

\bigskip
{\bf $\mathbf N$ is a complex space form and $\mathbf M$ is a Lagrangian
submanifold of $\mathbf N$\,:} Here we have
$\hol(N)=\fraku(T_oN)=\su(T_oN)\oplus\R\,\j$ and we aim to prove the equality
$\hol(TN|M)=\su(T_oN)$\,:

For this, notice that $\hol(N)_+$ resp.\ $\hol(N)_-$ is given by
$\su(T_oN)\cap\so(T_oN)_+$ resp.\ by $\R\,\j\oplus\su(T_oN)\cap\so(T_oN)_-$\,,
and hence the linear maps given by $\hol(N)_+\to\so(T_oM),\,A\mapsto
A|T_oM$ and $\hol(N)_+\to\so(\bot_oM),\,A\mapsto
A|\bot_oM$ both are isomorphisms. Therefrom, we easily verify that
$\frakc(\hol(N)_+)\cap\so(T_oN)_-=\R\,\j$ holds, hence we have
$\dim\big(\frakc(\hol(N)_+)\cap\hol(N)_-\big)=1$ and thus
$\hol(TN|M)_-=\{0\}$ is not possible by virtue of Part~(b) of Lemma~\ref{le:about_hol}.
Moreover,  $L$ is isomorphic to the
Grassmannian manifold of Lagrangian planes in $T_oN$\,, whose universal
covering space is a
product of an irreducible symmetric space and a 1-dimensional factor:
The corresponding decomposition of ``de\,Rham type'' is given by $\hol(N)_-=V_1\oplus
V_2$ with $V_1:=\R\,\j$ and $V_2:=\su(T_oN)\cap\so(T_oN)_-$\,. Thus Part~(a) of Lemma~\ref{le:about_hol} implies that $\hol(TN|M)_-$ is equal
to one of the spaces $\R\,\j$\,,
$\su(T_oN)_-$ or $\R\,\j\oplus\su(T_oN)_-$\,.
But $\j$ is orthogonal to $\hol(TN|M)$\,, as a consequence of Part~(c) of Lemma~\ref{le:about_hol}; therefore 
the only remaining possibility is $\hol(TN|M)=\su(T_oN)$\,.

\bigskip
{\bf $\mathbf N^{4n}$ is a quaternionic space form
(with $n\geq 2$) and $\mathbf M^{2n}$ is a totally complex submanifold of $\mathbf N$\,:} Here
the holonomy Lie algebra $\hol(N)$ is given by $\sp(T_oN)\oplus\rmQ$\,, where $\rmQ$
denotes the quaternionic structure at $o$\,. Choose a canonical basis
$\{\i,\j,\k\}$ with $\i(T_oM)=T_oM$\,, $\j(T_oM)=\bot_oM$ and $\k(T_oM)=\bot_oM$\,.
Let us see that that $\hol(TN|M)=\sp(n)\oplus\R\,\i$ holds:
 
We notice that $\hol(N)_+=\sp(T_oN)_+\oplus\R\,\i$
and $\hol(N)_-=\sp(T_oN)_-\oplus\R\,\j\oplus\R\,\k$ with $\sp(T_oN)_\pm:=\sp(T_oN)\cap\so(T_oN)_\pm$ and that the
linear maps $\sp(T_oN)_+\to\fraku(T_oM),\,A\mapsto
A|T_oM$ and $\sp(T_oN)_+\to\fraku(\bot_oM),\,A\mapsto
A|\bot_oM$ both are isomorphisms. Hence 
$\frakc(\sp(T_oN)_+\cap\so(T_oN)_-=\{\j,\k\}_\R$\,, thus
$\frakc(\hol(N)_+)\cap\so(T_oN)_-=\{0\}$\,, since $\i$ does not commute with
$\j$ or with $\k$\,. Therefore, $\hol(TN|M)_-=\{0\}$ again is not possible.
Moreover, $L$ is isomorphic to the
Grassmannian manifold of totally complex $2n$-planes in $T_oN$\,, whose
universal covering is a product of  two irreducible factors:
The corresponding decomposition of ``de\,Rham type'' is given by $\hol(N)_-=V_1\oplus
V_2$ with $V_1:=\sp(T_oN)_-$ and $V_2:=\{\j,\k\}_\R$\,.
Now the result follows by means of Part~(d) of Lemma~\ref{le:about_hol} combined with similar arguments as in the Lagrangian case.
\end{proof}

\subsection{\ldots in an ambient symmetric space of higher rank}
\label{s:higher_rank}
In this section, we will prove Theorem~\ref{th:hol_for_symmetric_M} in case the ambient space $N$ is of higher rank. 
As usual, let $o\in N$ be an origin, $K$ denote the isotropy subgroup of
$\Iso^0(N)$ at $o$\,, $\fraki(N)=\frakk\oplus\frakp$ the corresponding Cartan decomposition
and $\pi_2:\frakk\to\so(T_oN)$ the linearized isotropy representation (see Section~\ref{s:results}).

Suppose that $M\subset N$ is a full symmetric submanifold with $o\in M$\,. 
Then Theorem~\ref{th:Naitoh} ensures that there exists some $X\in\frakp$ with
$\ad(X)^3=\ad(X)$ (in the non-compact case) resp. $\ad(X)^3=-\ad(X)$ (in the
compact case) such that $M$ belongs to the family of symmetric submanifolds
associated with the symmetric R-space $\Ad(K)\,X$ according to Definition~\ref{de:BENT}\,; hence  $M=M_c$ for some $c\neq 0$ \,.
Then proposition~\ref{p:BENT} states that $T_oM=\pi_1(\frakp_\epsilon)$\,,
$\bot_oM=\pi_1(\frakp_0)$ and $\forall Y\in\frakp_\epsilon:\fetth(\pi_1(Y))=c\,\pi_2(J(Y))$
(where $J$ is the linear map defined in Lemma~\ref{le:ad(X)})\,.

\bigskip
\begin{lemma}\label{le:coincide}
By means of the identification $\frakk\cong\hol(N)$ from  Proposition~\ref{p:frakk_gleich_hol}, the
splittings
$\frakk=\frakk_0\oplus\frakk_\epsilon$ (see~\eqref{eq:decomposition_of_frakk_with_respect_to_ad(X)}) 
and $\hol(N)=\hol(N)_+\oplus\hol(N)_-$ (see~\eqref{eq:decomposition_of_hol(N)},\eqref{eq:hol(N)_pm})
are in correspondence with each other.
\end{lemma}
\begin{proof}
The result follows from the previous in combination with Equations~\eqref{eq:ad_frakp},~\eqref{eq:hol(N)_pm} and~\eqref{eq:correspondence}.
\end{proof} 

\bigskip
\begin{proposition}\label{p:about_hol_of_N}
We always have $\fetth(T_oM)=\hol(TN|M)_-$\,.
Furthermore, the decomposition $\hol(TN|M)=\hol(TN|M)_+\oplus\hol(TN|M)_-$
mentioned in Proposition~\ref{p:frakk_bot} is given by $\hol(TN|M)=\hol(N)_+\oplus[\hol(N)_-,\hol(N)_+]$\,.
Therefore, $\hol(TN|M)$ contains the ideal $[\hol(N),\hol(N)]$\,.
\end{proposition}
\begin{proof}
Because of~\eqref{eq:ad(X)_2} combined with Proposition~\ref{p:BENT}, we have
$\fetth(T_oM)=\pi_2(\frakk_\epsilon)=\hol(N)_-$\,, where the second equality
follows from Lemma~\ref{le:coincide}; therefore, as a consequence
of~\eqref{eq:hol_minus} and Proposition~\ref{p:frakk_bot}, we have $\hol(TN|M)=\hol(N)_+\oplus[\hol(N)_-,\hol(N)_+]$\,.
The last assertion is now seen from 
\begin{align*}
[\hol(N),\hol(N)]&=[\hol(N)_+,\hol(N)_+]+[\hol(N)_-,\hol(N)_-]+[\hol(N)_-,\hol(N)_+]\\
&\subset\hol(N)_+\oplus[\hol(N)_-,\hol(N)_+]\;.
\end{align*}
\end{proof}

\bigskip
Suppose that $N$ is a Hermitian symmetric space. 
Then we have the splitting $\frakk=\frakc\oplus [\frakk,\frakk]$ with a
one-dimensional factor $\frakc$\,, the center of $\frakk$\,. Moreover, there exists
$Z\in\frakc$ with $\pi_2(Z)=\j$\,, where the latter denotes the complex
structure of $T_oN$ (cf.~\cite{BCO},~A.\,4).

\bigskip
\begin{lemma}\label{le:frakc}
We have $\frakc\subset\frakk_\epsilon$\,. Furthermore, $M$ is a Lagrangian submanifold of $N$\,.
\end{lemma}
\begin{proof}
Let $\sigma_o^\bot$ denote the extrinsic symmetry of $M$ at $o$ and 
consider the (second) involution $\tau:=\Ad(\sigma_o^\bot)$ on $\fraki(N)$\,.
Then we have $\tau(\frakk)=\frakk$ and (in the notation of Section~\ref{s:extrinsically_hol}) 
$$
\forall X\in\frakk:\pi_2(\tau\,X)=\sigma^\bot\circ \pi_2(X)\circ\sigma^\bot\;.
$$
Thus Lemma~\ref{le:coincide} implies that the splitting $\frakk=\frakk_0\oplus\frakk_\epsilon$ is also the
decomposition of $\frakk$ into the $+1$ and $-1$ eigenspaces of
$\tau|\frakk$\,; therefore, and since $\tau|\frakk$ is a Lie algebra involution, it maps 
$\frakc$ onto itself, hence either $\frakc\subset\frakk_0$ or $\frakc\subset\frakk_\epsilon$\,.   
By contradiction, if we had $\frakc\subset\frakk_0$\,, then $[X,Z]=0$ according 
to~\eqref{eq:decomposition_of_frakk_with_respect_to_ad(X)}, thus
$$\j(\pi_1(X))=\pi_2(Z)\,\pi_1(X)\stackrel{\eqref{eq:correspondence}}=\pi_1(\ad(Z)\,X)=0$$
(because $Z$ belongs to the center of $\frakk$);
therefore, and since $\j$ is the complex structure of
$T_oN$\,, we have $\pi_1(X)=0$ and hence $X=0$\,, which is not possible.
Thus we obtain $\frakc\subset\frakk_\epsilon$ and therefore
$\j\in\pi_2(\frakk_\epsilon)\subset\so(T_oN)_-$\,, which, by virtue of~\eqref{eq:odd},
implies that $\j$ maps $T_oM$ to $\bot_oM$ and vice versa. Since $M$ is an extrinsically homogeneous submanifold of $N$\,,
we hence see that $M$ is already a Lagrangian submanifold.
\end{proof}

\begin{proof}[Proof of Theorem~\ref{th:hol_for_symmetric_M} in case $N$ is of higher rank] If $N$ is of Hermitian type, then $M$ is a Lagrangian submanifold of $N$ according to Lemma~\ref{le:frakc}.
Furthermore, Proposition~\ref{p:about_hol_of_N} implies that
$[\hol(N),\hol(N)]\subset\hol(TN|M)\subset\hol(N)=[\hol(N),\hol(N)]\oplus\R\,\j$\,. Because the complex structure $\j$ is orthogonal to $\hol(TN|M)$ as a result of
Lemma~\ref{le:about_hol}, in fact we have $\hol(TN|M)=[\hol(N),\hol(N)]$\,.
Moreover, $\j=\pi_2(J)\in\pi_2(\frakk)_-=\fetth(T_oM)$ by means of
Lemma~\ref{le:frakc} combined with Proposition~\ref{p:about_hol_of_N}.

If $N$ is not of Hermitian type, then $\frakk$ is semisimple (because $N$ is not of Hermitian type and
hence the center of $\fraki(N)$ is trivial) and thus we have
$\frakk=[\frakk,\frakk]$\,. Therefore, $\hol(TN|M)=\hol(N)$ follows by virtue of Proposition~\ref{p:about_hol_of_N}.
\end{proof}

\vspace{2cm}
\begin{center}

 \qquad
 \parbox{60mm}{Tillmann Jentsch\\
  Mathematisches Institut\\
  Universit{\"a}t zu K{\"o}ln\\
  Weyertal 86-90\\
  D-50931 K{\"o}ln, Germany\\[1mm]
  \texttt{tjentsch@math.uni-koeln.de}}

\end{center}
\end{document}